\documentclass[12pt]{article}
\usepackage{graphicx}
\usepackage{amsmath,amsthm,amssymb,enumerate,esint}
\usepackage{color}
\usepackage[left=2cm,right=2cm,top=3.5cm,bottom=3.5cm]{geometry}
\usepackage{color}
\usepackage{hyperref}
\hypersetup{
  colorlinks = true, 
  urlcolor = blue, 
  linkcolor = blue, 
  citecolor = blue 
}
\allowdisplaybreaks

\newtheorem{Theorem}{Theorem}[section]
\newtheorem{Proposition}[Theorem]{Proposition}
\newtheorem{Lemma}[Theorem]{Lemma}
\newtheorem{Corollary}[Theorem]{Corollary}

\theoremstyle{definition}
\newtheorem{Definition}[Theorem]{Definition}

\newtheorem{Remark}[Theorem]{Remark}

\newcommand{\bTheorem}[1]{
	\begin{Theorem} \label{T#1} }
	\newcommand{\eT}{\end{Theorem}}

\newcommand{\bProposition}[1]{
	\begin{Proposition} \label{P#1}}
	\newcommand{\eP}{\end{Proposition}}

\newcommand{\bLemma}[1]{
	\begin{Lemma} \label{L#1} }
	\newcommand{\eL}{\end{Lemma}}

\newcommand{\bCorollary}[1]{
	\begin{Corollary} \label{C#1} }
	\newcommand{\eC}{\end{Corollary}}

\newcommand{\bRemark}[1]{
	\begin{Remark} \label{R#1} }
	\newcommand{\eR}{\end{Remark}}
	
\newcommand{\bDefinition}[1]{
	\begin{Definition} \label{D#1} }
	\newcommand{\eD}{\end{Definition}}
	
\newcommand{\bFormula}[1]{
	\begin{equation} \label{#1}}
	\newcommand{\eF}{\end{equation}}

\newcommand{\Ma}{{\rm Ma}}
\newcommand{\Fr}{{\rm Fr}}
\newcommand{\Al}{{\rm Al}}
\newcommand{\Ov}[1]{\overline{#1}}
\newcommand{\vc}[1]{{\bf #1}}
\newcommand{\intO}[1]{\int_{\Omega} #1 \ \dx}
\renewcommand{\S}{\mathbb{S}}
\newcommand{\eps}{\varepsilon}
\newcommand{\del}{\partial}
\newcommand{\vB}{\vc{B}}

\newcommand{\vr}{\varrho}
\newcommand{\vu}{\vc{u}}
\newcommand{\vm}{\vc{m}}

\newcommand{\vt}{\vartheta}

\newcommand{\vtB}{\vt_B}
\newcommand{\vU}{\vc{U}}
\newcommand{\dd}{\ \mathrm{d}}
\newcommand{\Div}{{\rm div}_x}

\newcommand{\Grad}{\nabla_x}

\newcommand{\Curl}{{\bf curl}_x}
\newcommand{\dx}{\,{\rm d} {x}}
\newcommand{\dt}{\,{\rm d} t }
\newcommand{\R}{\mathbb{R}}
\newcommand{\I}{\mathbb{I}}
\newcommand{\br}{ \nonumber \\ }
\newcommand{\weak}{\rightharpoonup}
\newcommand{\weakstar}{\stackrel{\ast}{\rightharpoonup}} 

\begin{document}

\title{Rigorous derivation of magneto-Oberbeck-Boussinesq approximation with non-local temperature term}

\author{Piotr Gwiazda
		\thanks{The work of P. G. was partially supported by  National Science Centre
			(Poland),  agreement no 2021/43/B/ST1/02851.}
		\and Florian Oschmann
			\thanks{The work of F. O. has been supported by the Czech Science Foundation (GA\v CR) project 22-01591S, and the Czech Academy of Sciences project L100192351. The Institute of Mathematics, CAS is supported by RVO:67985840.}
		\and Aneta Wr\'oblewska-Kami\'nska
			\thanks{The work of A. W.-K. was supported by National Science Centre Poland, grant no. UMO-2020/38/E/ST1/00469}			
			}

\date{}

\maketitle

\medskip

\centerline{Institute of Mathematics of Polish Academy of Sciences}
\centerline{\'Sniadeckich 8, 00-956 Warszawa, Poland}

\medskip

\centerline{Institute of Mathematics of the Academy of Sciences of the Czech Republic}
\centerline{\v Zitn\' a 25, CZ-115 67 Praha 1, Czech Republic}

\begin{abstract}
We consider a general compressible, viscous, heat and magnetically conducting fluid described by the compressible Naiver–Stokes–Fourier system coupled with induction equation. In particular, we do not assume conservative boundary conditions for the temperature and allow heating or cooling on the surface of the domain. We are interested in the mathematical analysis when the Mach, Froude, and Alfv\'en numbers are small, converging to zero at a specific rate. We give a rigorous mathematical justification that in the limit, in case of low stratification, one obtains a modified Oberbeck-Boussinesq-MHD system with a non-local term or a non-local boundary condition for the temperature deviation. Choosing a domain confined between parallel plates, one finds also that the flow is horizontal, and the magnetic field is perpendicular to it. The proof is based on the analysis of weak solutions to a primitive system and the relative entropy method.
\end{abstract}


\tableofcontents

\section{Introduction}
Let $\Omega$ be the three-dimensional periodic strip
\begin{align*}
\Omega = \mathbb{T}^2 \times (0,1), \ \mathbb{T}^2 = \left( [-1,1] \Big|_{\{-1,1\}} \right)^2.
\end{align*}
Let $T>0$, and let $(0,T) \times \Omega$ be a time-space cylinder on which we study the \emph{compressible magneto-hydrodynamic (MHD) system} describing the behavior of four quantities: density of a fluid $\vr: (0,T)\times \Omega \to \R$, velocity field $\vu: (0,T)\times \Omega \to \R^3$, absolute temperature $\vt: (0,T) \times \Omega \to \R$, magnetic field $\vc{B}: (0,T) \times \Omega \to \R^3$, and consisting of equation of continuity, momentum equation, induction equation, and  entropy balance, respectively: 
%
%
	\begin{equation} \label{p1}
		\partial_t \vr + \Div (\vr \vu) = 0,
	\end{equation}
%
%
	\begin{equation} \label{p2}
		\partial_t (\vr \vu) + \Div (\vr \vu \otimes \vu) - \Div \mathbb{S}(\vt, \Grad \vu) + \frac{1}{\Ma^2} \Grad p (\vr, \vt) = \frac{1}{\Fr^2} \vr \Grad G + \frac{1}{\Al^2} \Curl  \vc{B} \times \vc{B},
	\end{equation}
%
%
\begin{equation} \label{p3}
	\partial_t \vc{B} + \Curl (\vc{B} \times \vu ) + \Curl (\zeta (\vt) \Curl  \vc{B} ) = 0,\
	\Div \vc{B} = 0,
	\end{equation}
%
	\begin{align}
	\partial_t (\vr s(\vr, \vt)) &+ \Div (\vr s(\vr, \vt) \vu) + \Div \left( \frac{\vc{q}(\vt, \Grad \vt)}{\vt} \right) \br &=
	\frac{1}{\vt} \left( \Ma^2 \mathbb{S}(\vt, \Grad \vu) : \Grad \vu - \frac{\vc{q}(\vt, \Grad \vt) \cdot \Grad \vt}{\vt} + \frac{\Ma^2}{\Al^2} \zeta(\vt) |\Curl \vB |^2 \right). \label{p4}
\end{align}

Here, the pressure $p = p(\vr, \vt)$, the entropy $s = s(\vr, \vt)$, and the internal energy $e=e(\vr, \vt)$ are interrelated through \emph{Gibbs' equation}
\begin{equation} \label{w1}
	\vt D s = De + p D \left( \frac{1}{\vr} \right).
\end{equation}
Moreover, the numbers $\Ma = \Ma(\eps), \, \Fr = \Fr(\eps), \, \Al = \Al(\eps)$ are the Mach, Froude, and Alfv\'en number, respectively, which we assume to depend on a small parameter $\eps>0$. More precisely, we consider the case of \emph{low stratification} in the sense that
\begin{align*}
\Ma=\Al=\eps, \ \Fr=\sqrt{\eps}.
\end{align*}

The right-hand side of equation \eqref{p4} represents the entropy production rate, which by the Second law of thermodynamics must be non-negative. Thus, we consider a \emph{Newtonian fluid}, 
with the viscous stress tensor given by
\begin{equation} \label{c7}
	\mathbb{S}(\vt, \Grad \vu) = \mu (\vt) \left( \Grad \vu + \Grad^T \vu - \frac23 \Div \vu \mathbb{I} \right) +
	\eta(\vt) \Div \vu \mathbb{I},
\end{equation}
where the viscosity coefficients $\mu > 0$ and $\eta \geq 0$ are continuously differentiable functions of the temperature. Similarly, the heat flux obeys \emph{Fourier's law},
\begin{equation} \label{c9}
	\vc{q}(\vt, \Grad \vt)= - \kappa (\vt) \Grad \vt,
\end{equation}	
where the heat conductivity coefficient $\kappa > 0$ is a continuously differentiable function of the temperature. Lastly, $\zeta \geq 0$ is the magnetic diffusivity.\\

We supplement the MHD system \eqref{p1}--\eqref{p4} with slip boundary conditions for the velocity and inhomogeneous Dirichlet boundary condition for temperature, namely 
\begin{equation} \label{p7}
	\vu \cdot \vc n|_{\partial \Omega} = 0, \quad (\S(\vt, \Grad \vu) \cdot \vc n) \times \vc n|_{\del \Omega} = 0, \quad \vt|_{\partial \Omega} = \Ov{\vt} + \eps \vtB, \quad \Ov{\vt} \in (0,\infty), \ \vt_B \in C^2(\del \Omega).
\end{equation}
For the magnetic field we assume a perfectly isolating boundary, meaning
\begin{equation}\label{p72}
\vB \times \vc n|_{\del \Omega} = 0.
\end{equation}

Finally, the initial conditions read
\begin{equation}\label{p73}
\vr(0, \cdot) = \vr_0 \in L^\frac{5}{3}(\Omega), \ (\vr \vu)(0, \cdot) = \vc m_0 \in L^\frac{5}{4}(\Omega), \ \vt(0, \cdot)=\vt_0 \in W^{1,2}(\Omega), \ \vB(0, \cdot) = \vB_0 \in L^2(\Omega).
\end{equation}

The existence of weak solutions to system \eqref{p1}--\eqref{p4} and weak-strong-uniqueness property was shown in \cite{FGKS2023}, see also the paper \cite{DucometFeireisl2006}, where the steps for existence proof are carried out, and we can do similarly by taking the ballistic energy to incorporate the temperature Dirichlet boundary conditions as introduced in \cite{ChaudhuriFeireisl2022}.\\

The same system as considered here without gravitational force ($G=0$) and with no-flux conditions for the temperature ($\vc q \cdot \vc n = 0$) has been considered in a bounded domain $\Omega \subset \R^3$ by Kuku\v cka in \cite{Kukuvcka2011}. The author's limiting system is similar to the one we will derive in the next section, but without non-local term. Furthermore, the incompressible Euler system without heat equation has been considered in \cite{JiangJu2019} for the torus case $\Omega = \mathbb{T}^2 \times (0,1)$. There, the outcome is a horizontal Euler flow, where $\vU=(U_1, U_2, 0)(t, x_1, x_2)$ and $\vB^1 = (0,0,b^1)(t, x_1, x_2)$ are the limiting functions. We will see that this structure holds in our case as well.\\

Besides the mentioned mathematical literature, there are several works from the physics community. Without being exhaustive, let us mention the work \cite{SpiegelWeiss1982}, who derived the magneto-Oberbeck-Boussinesq approximation rather formally, but also investigates unstabilities, magnetoconvection, and magnetic buoyancy. Subsequent results are given by Corfield \cite{Corfield1984}, giving a more refined derivation by splitting in ``parallel'' and ``orthogonal'' directions. Force-free equations, where the Lorentz force $\Curl \vB \times \vB = 0$, are studied in \cite{Chandrasekhar1956, ChandrasekharWoltjer1958, Tassi2008}. For an introduction to the derivation of the Oberbeck-Boussinesq approximation without magnetic field, we refer to \cite{Mizerski2021} and the references therein.



\section{Formal limit}\label{sec:formLim}
Let us formally derive the limit system of the compressible MHD system \eqref{p1}-\eqref{p4} as $\eps \to 0$. To start we may write a formal expansion 
\begin{equation*}
\begin{split}
\vr & = \Ov{\vr}  + \eps \vr^{1} + \dots\\
\vu & = \vc{U}  + \eps \vu^{1} + \dots\\
\vt & = \Ov{\vt}  + \eps \vt^{1} + \dots\\
\vB & = \Ov{\vB}  + \eps \vB^{1} + \dots\\
\end{split}
\end{equation*}
Inserting the above expansion to the system \eqref{p1}-\eqref{p4} and regrouping terms with respect to powers of $\eps$ we obtain that $\Grad p(\Ov{\vr}, \Ov{\vt}) = \Curl \Ov{\vB} \times \Ov{\vB}$.
Since we are interested in global solutions, the necessary estimate needs to be obtained from the dissipation equation and in particular the entropy production rate needs to be kept small of order $\eps^2$. Consequently, the quantities $\Grad \vt$ and $\Curl \vB$ vanish in the asymptotic limit $\eps \to 0$. It is natural then to assume that $\Ov\vt$ is a positive constant. As $\Div \vB = 0$, we may expect that  $\Ov{\vB} = \Ov{\vB}(t) \in \R^3$; the induction equation then shows that $\Ov{\vB}$ is indeed independent of time and hence an overall constant. Imposing further the boundary conditions for $\vB$, we infer
\begin{align*}
    \Ov{\vB} = (0, 0, \Ov{b})^T \ \text{for some} \ \Ov{b} \in \R.
\end{align*}
%
This also enforces that $\Ov{\vr} > 0$ will be a constant.

Next, the zero-th order terms in the momentum equation \eqref{p2} give rise to
\begin{align}\label{firstGrad}
\Grad p(\vr, \vt) = \eps \vr \Grad G + \Curl \vB \times \vB + O(\eps^2).
\end{align}
Accordingly, a solution to the MHD system can be written as $$\vr = \Ov{\vr} + \eps \vr_\eps^1, \quad \vt = \Ov{\vt} + \eps \vt_\eps^1,\quad  \vB = \Ov{\vB} + \eps \vB_\eps^1,$$ where $\Ov{\vr}, \Ov{\vt}>0$ are positive constants, and $\Ov{\vB} = (0,0, \Ov{b})^T$ for some $\Ov{b} \in \R$.

Anticipating now the convergences $(\vr_\eps^1, \vt_\eps^1, \vB_\eps^1) \to (\vr^1, \vt^1, \vB^1)$ (in some sense), we find
\begin{align}\label{OB1}
\del_\vr p(\Ov{\vr}, \Ov{\vt}) \Grad \vr^1 + \del_\vt p(\Ov{\vr}, \Ov{\vt}) \Grad \vt^1 = \Ov{\vr}\Grad G + \Curl \vB^1 \times \Ov{\vB}.
\end{align}
For structural reasons, we see that there must exist some function $A$ such that $- \Grad A = \Curl \vB^1 \times \Ov{\vB}$, and without loss of generality, we might assume $\intO{A} = 0$.
Removing now gradients in \eqref{OB1}, this yields
\begin{align*}
\del_\vr p(\Ov{\vr}, \Ov{\vt}) \vr^1 + \del_\vt p(\Ov{\vr}, \Ov{\vt}) \vt^1 + A = \Ov{\vr}G + \chi(t)
\end{align*}
for some spatially homogeneous function $\chi$. By conservation of mass, we shall suppose $\int_\Omega \vr^1 \dd x = 0$, and without loss of generality, we can also assume $\int_\Omega G \dd x = 0$. In turn, we find
\begin{equation}\label{chi_t}
\chi(t) = \del_\vt p(\Ov{\vr}, \Ov{\vt}) \fint_\Omega \vt^1 \dd x,
\end{equation}
thus the \emph{magnetic Boussinesq relation} reads
\begin{align}\label{OB}
\del_\vr p(\Ov{\vr}, \Ov{\vt}) \vr^1 + \del_\vt p(\Ov{\vr}, \Ov{\vt}) \vt^1 + A = \Ov{\vr}G + \del_\vt p(\Ov{\vr}, \Ov{\vt}) \fint_\Omega \vt^1 \dd x.
\end{align}

Focusing on the induction equation \eqref{p3}, we get
\begin{align}\label{mag}
\eps \del_t \vB^1 + \Curl(\Ov{\vB} \times \vU) + \eps \Curl(\vB^1 \times \vU) + \eps \Curl(\zeta(\Ov{\vt} + \eps \vt^1) \Curl \vB^1) = 0.
\end{align}

Collecting finally terms in \eqref{mag} for $\eps \to 0$, we infer $\Curl(\Ov{\vB} \times \vU) = 0$, and the induction equation for $\vB^1$ reads
\begin{align*}
\del_t \vB^1 + \Curl(\vB^1 \times \vU) + \Curl (\zeta(\Ov{\vt}) \Curl \vB^1) = 0.
\end{align*}

Note that this implies by $\Ov{\vB} = (0, 0, \Ov{b})^T$ and $\Div \vU = 0$ that
\begin{align*}
    0 = \Curl(\Ov{\vB} \times \vU) = (\Ov{\vB} \cdot \Grad) \vU = \Ov{b} \del_3 \vU,
\end{align*}
hence $\vU$ is independent of the third coordinate, which together with the boundary condition $\vU \cdot \vc n |_{\del \Omega} = 0$ implies
\begin{align*}
    \vU = (U_1, U_2, 0)^T(t, x_1, x_2).
\end{align*}

As for the continuity equation, we infer
\begin{align*}
    0 = \Ov{\vr} \Div \vU + \eps (\del_t \vr^1 + \Div(\vr^1 \vU)),
\end{align*}
hence $\Div \vU = 0$ as expected, and also $\del_t \vr^1 + \Div(\vr^1 \vU) = 0$.\\

The momentum equation is derived as follows: first, we have from \eqref{p2}
\begin{align*}
    \vr_\eps (\del_t \vu_\eps + \vu_\eps \cdot \Grad \vu_\eps) - \Div \S(\vt_\eps, \Grad \vu_\eps) = \frac{1}{\eps}\vr_\eps \Grad G - \frac{1}{\eps^2} \Grad p(\vr_\eps, \vt_\eps) + \frac{1}{\eps^2} \Curl \vB_\eps \times \vB_\eps.
\end{align*}
The left hand-side of the above equation tends to
\begin{align*}
    \Ov{\vr} (\del_t \vU + \vU \cdot \Grad \vU) - \Div \S(\Ov{\vt}, \Grad \vU).
\end{align*}
For the right hand-side, we write using the expansions of $\vr_\eps, \vt_\eps, \vu_\eps$, and $\vB_\eps$
\begin{align*}
    &\frac{1}{\eps}\vr_\eps \Grad G - \frac{1}{\eps^2} \Grad p(\vr_\eps, \vt_\eps) + \frac{1}{\eps^2} \Curl \vB_\eps \times \vB_\eps \br 
    &= \frac{1}{\eps} (\Ov{\vr} \Grad G - \del_\vr p(\Ov{\vr}, \Ov{\vt}) \Grad \vr^1 - \del_\vt p(\Ov{\vr}, \Ov{\vt}) \Grad \vt^1 + \Curl \vB^1 \times \Ov{\vB}) \br 
    &\qquad - \Grad  \Pi + \vr^1 \Grad G + \Curl \vB^1 \times \vB^1 + O(\eps) \br 
    &= \vr^1 \Grad G - \Grad \Pi + \Curl \vB^1 \times \vB^1 + O(\eps),
\end{align*}
where we used that the term in brackets is zero due to $\Curl \vB^1 \times \Ov{\vB} = - \Grad A$ and Boussinesq relation \eqref{OB}, and the term $\Grad \Pi$ corresponds to the $\eps^0$-order term of the expansion of $\eps^{-2} \Grad p(\vr_\eps, \vt_\eps)$. Hence, the limiting momentum equation reads
\begin{align*}
    \Ov{\vr} (\del_t \vU + \vU \cdot \Grad \vU) - \Div \S(\Ov{\vt}, \Grad \vU) + \Grad \Pi = \vr^1 \Grad G + \Curl \vB^1 \times \vB^1.
\end{align*}

To handle the entropy equation \eqref{p4}, we insert our expansion into the entropy equation, use that the entropy production rate is small of order $O(\eps^2)$, and perform the limit $\eps \to 0$ to obtain
\begin{align*}
    \Ov{\vr} \big( \del_\vr s(\Ov{\vr}, \Ov{\vt}) \del_t \vr^1 + \del_\vt s(\Ov{\vr}, \Ov{\vt}) \del_t \vt^1 + \vU \cdot \Grad ( \del_\vr s(\Ov{\vr}, \Ov{\vt}) \vr^1 + \del_\vt s(\Ov{\vr}, \Ov{\vt}) \vt^1 ) \big) - \frac{\kappa(\Ov{\vt})}{\Ov{\vt}} \Delta_x \vt^1 = 0.
\end{align*}
Finally, we use \eqref{OB} to express the density $\vr^1$ by the corresponding other terms and replace it in the above equation, and use also that by Gibbs' relation, one finds $\vt \del_\vt s(\vr, \vt) = \del_\vt e(\vr, \vt)$ and $\vt \del_\vr s(\vr, \vt) = -\vt \vr^{-2} \del_\vt p(\vr, \vt)$. Similar to the case without magnetic field (see \cite{BellaFeireislOschmann2023a}), the function $\chi$ in \eqref{OB} shall give rise to non-local effects. As a matter of fact, we will show that the limiting system will be the \emph{modified Oberbeck-Boussinesq MHD (OBM) system}
\begin{align}\label{target}
\Div \vU = \Div \vB^1 &= 0,\br
\del_t \vr^1 + \Div(\vr^1 \vU) &= 0, \br 
\Ov{\vr} \big( \del_t \vU + \vU \cdot \Grad \vU \big) - \Div \S(\Ov{\vt}, \Grad \vU) + \Grad \Pi &= \vr^1 \Grad G + \Curl \vB^1 \times \vB^1, \br
\Ov{\vr}c_p \big( \del_t \vt^1 + \vU \cdot \Grad \vt^1 \big) - \Ov{\vr} \Ov{\vt} \alpha \vU \cdot \Grad G - \Ov{\vt} \alpha (\del_t A + \vU \cdot \Grad A) &= \kappa(\Ov{\vt})\Delta_x \vt^1 + \Ov{\vt} \alpha \del_\vt p(\Ov{\vr}, \Ov{\vt}) \del_t \fint_\Omega \vt^1 \dd x,\br
\del_t \vB^1 + \Curl( \vB^1 \times \vU) + \Curl( \zeta(\Ov{\vt}) \Curl \vB^1) &= 0,
\end{align}
together with the Boussinesq relation \eqref{OB}, and the classical definitions of the coefficient of thermal expansion
and the specific heat at constant pressure, respectively,
\begin{align}\label{Coeff} 
	\alpha = \alpha(\Ov{\vr}, \Ov{\vt} ) \equiv \frac{1}{\Ov{\vr}}  \frac{\partial_\vt p(\Ov{\vr}, \Ov{\vt} ) }{\partial_\vr p(\Ov{\vr}, \Ov{\vt} )},  && 
	c_p = c_p (\Ov{\vr}, \Ov{\vt} ) \equiv \partial_\vt e(\Ov{\vr}, \Ov{\vt} ) + \frac{\Ov{\vt} \alpha}{\Ov{\vr}} \partial_\vt p(\Ov{\vr}, \Ov{\vt} ).
\end{align}

We complete system \eqref{target} with the corresponding boundary conditions
\begin{align}\label{targetBC}
\vU \cdot \vc n|_{\del \Omega} = 0, \ (\S(\Ov{\vt}, \Grad \vU) \cdot \vc n)\times \vc n|_{\del \Omega} = 0, \ \vt^1|_{\del \Omega} = \vt_B, \ \vB^1 \times \vc n|_{\del \Omega} = 0.
\end{align}

By the form of $\vU$, we see that the boundary conditions for the velocity are fulfilled trivially.

\subsection{On the structure of magnetic field and heat equation}
Let us have a closer look on $\vB^1$. First, as $\Curl \vB^1 \times \Ov{\vB} = - \Grad A$, we have
\begin{align*}
    0 = - \Curl \Grad A = \Curl (\Curl \vB^1 \times \Ov{\vB}) = (\Ov{\vB} \cdot \Grad) \Curl \vB^1 = \Ov{b} \del_3 \Curl \vB^1
\end{align*}
as both $\Ov{\vB}$ and $\Curl \vB^1$ are divergence-free, and $\Ov{\vB}$ is constant. Moreover, as $\Div \vB^1 = 0$, we infer $\Div \del_3 \vB^1 = 0$. Next, by $\Curl \del_3 \vB^1 = 0$ and the boundary condition $\vB^1 \times \vc n = 0$ on $\del \Omega$, it is reasonable to assume that also $(\del_3 \vB^1) \times \vc n = 0$ there. In turn, Friedrichs' inequality (see, e.g., \cite[Theorem~3.1]{vonWahl1992}) gives $\Grad \del_3 \vB^1 = 0$ in $\Omega$ and hence $\vB^1 = \vc c_1 x_3 + \vc c_2(x_1, x_2)$ for some (spatially) constant vector $\vc c_1 \in \R^3$. The boundary conditions (used separately on $x_3=0$ and $x_3=1$, respectively) then give $c_{2,1} = c_{2,2} = 0$ and $c_{1,1} + c_{2,1} = c_{1,2} + c_{2,2} = 0$, hence $\vB^1 = (0, 0, c_1 x_3 + c_2(x_1, x_2))^T$. Solenoidality now enforces $c_1 = 0$ and thus
\begin{align*}
    \vB^1 = (0, 0, b^1)^T(t, x_1, x_2).
\end{align*}
Note that this means that the only boundary condition surviving in \eqref{targetBC} is the one for the temperature. Moreover, in this case
\begin{align*}
    \Curl(\vB^1 \times \vU) = \Div(b^1 \vU) \vc e_3, && \Curl (\zeta(\Ov{\vt}) \Curl \vB^1) = - \zeta(\Ov{\vt}) \Delta_x b^1 \vc e_3,
\end{align*}
such that the induction equation takes the scalar form
\begin{align*}
    \del_t b^1 + \Div(b^1 \vU) - \zeta(\Ov{\vt}) \Delta_x b^1 = 0.
\end{align*}

\begin{Remark}\label{rem:formBheat}
In turn, we may write $\Curl \vB^1 \times \Ov{\vB} = - \Grad(\Ov{\vB} \cdot \vB^1)$, so $A = \Ov{\vB} \cdot \vB^1$, meaning also that the magnetic tension $(\Ov{\vB} \cdot \Grad) \vB^1 = 0$, which is meaningful as the magnetic field lines are straight. Moreover, then $\Curl \vB^1 \times \vB^1 = - \Grad \frac12 |\vB^1|^2$ which can be absorbed into the pressure term $\Grad \Pi$ and correctly gives the total pressure as sum of gas and magnetic pressure.
\end{Remark}

As for the heat equation, we see now that $A = \Ov{\vB} \cdot \vB^1 = \Ov{b} b^1$ fulfills the induction equation and we may write in more compact form
\begin{align*}
&\kappa(\Ov{\vt})\Delta_x \vt^1 + \Ov{\vt} \alpha \del_\vt p(\Ov{\vr}, \Ov{\vt}) \del_t \fint_\Omega \vt^1 \dd x \br 
    &=\Ov{\vr}c_p \big( \del_t \vt^1 + \vU \cdot \Grad \vt^1 \big) - \Ov{\vr} \Ov{\vt} \alpha \vU \cdot \Grad G - \Ov{\vt} \alpha (\del_t A + \vU \cdot \Grad A) \br  
    &=\Ov{\vr}c_p \big( \del_t \vt^1 + \vU \cdot \Grad \vt^1 \big) - \Ov{\vr} \Ov{\vt} \alpha \vU \cdot \Grad G - \Ov{\vt} \alpha \zeta(\Ov{\vt}) \Delta_x A.
\end{align*}
This is the heat equation occurring in \cite{SpiegelWeiss1982} (without non-local temperature term).




\section{Weak formulation}
\label{w}

We start by a classical list of structural hypotheses imposed on the equation of state and the transport coefficients.

\subsection{Equation of state}

Our choice of the equation of state is motivated by \cite[Chapter 4]{FeireislNovotny2022}, and we take the same assumptions as made in \cite{FGKS2023}. In particular, we suppose that
\begin{align}
	p(\vr, \vt) &= p_M (\vr, \vt) + p_R (\vt), \ p_M(\vr,\vt) = \vt^{\frac{5}{2}} P \left( \frac{\vr}{\vt^{\frac{3}{2}}  } \right), \ p_R (\vt) = \frac{a}{3} \vt^4,\br
	e(\vr, \vt) &= e_M (\vr, \vt) + e_R (\vr, \vt),\ \vr e_M(\vr,\vt) =  \frac{3}{2} \vt^{\frac{5}{2}} P \left( \frac{\vr}{\vt^{\frac{3}{2}}  } \right),\ \vr e_R (\vr ,\vt) = a \vt^4, \ a > 0,
	\label{c1}	
\end{align}
where the function $P \in C^1[0,\infty)$ satisfies
\begin{equation} \label{c2}
	P(0) = 0,\ P'(Z) > 0 \ \mbox{for}\ Z \geq 0,\ 0 < \frac{ \frac{5}{3} P(Z) - P'(Z) Z }{Z} \leq c \ \mbox{for}\ Z > 0.
\end{equation} 	
This implies that $Z \mapsto P(Z)/ Z^\frac{5}{3}$ is decreasing, and we suppose
\begin{equation} \label{c3}
	\lim_{Z \to \infty} \frac{ P(Z) }{Z^\frac{5}{3}} = p_\infty > 0.
\end{equation}
Note that $p_M$, $e_M$ satisfy the relation
\[
p_M = \frac{2}{3} \vr e_M,
\]
$p_R$ is the radiation pressure with the associated internal energy $e_R$, and \eqref{c3} corresponds to the
presence of electron pressure under the degenerate gas regime.
It follows from \eqref{c2} that $p$ and $e$ satisfy the \emph{hypothesis of thermodynamic stability}:
\begin{equation} \label{c4}
	\del_\vr p (\vr, \vt) > 0, \qquad \del_\vt e (\vr, \vt) > 0.
\end{equation}

In accordance with \eqref{w1}, the entropy takes the form
\begin{equation} \label{c5}
	s(\vr, \vt) = s_M(\vr, \vt) + s_R (\vr, \vt),\ s_M(\vr, \vt) = \mathcal{S} \left( \frac{\vr}{\vt^{\frac{3}{2}} } \right),\ \vr s_R(\vr, \vt) = \frac{4a}{3} \vt^3,
\end{equation}
where
\begin{equation} \label{c6}
	\mathcal{S}'(Z) = -\frac{3}{2} \frac{ \frac{5}{3} P(Z) - P'(Z) Z }{Z^2} < 0.
\end{equation}

\begin{Remark}
    For \eqref{c1}, the constant $a>0$ is crucial in the existence theory, as it rules out possible oscillations of the temperature in the vacuum zone $\{\vr = 0\}$. As it does not influence the final form of the system, we might instead consider a function $a=a(\eps) \to 0$ as $\eps \to 0$.
\end{Remark}

\subsection{Transport coefficients}

We suppose the viscosity coefficients in the Newtonian stress $\mathbb{S}(\vt, \Grad \vu)$ are continuously
differentiable functions of the temperature satisfying
\begin{align}
	0 < \underline{\mu} \left(1 + \vt \right) &\leq \mu(\vt) \leq \Ov{\mu} \left( 1 + \vt \right),\
	|\mu'(\vt)| \leq c \ \mbox{for all}\ \vt \geq 0, \br
	0 &\leq  \eta(\vt) \leq \Ov{\eta} \left( 1 + \vt \right).
	\label{c8}
\end{align}

Similarly, the heat conductivity coefficient in the Fourier heat flux $\vc{q}(\vt, \Grad \vt) = - \kappa(\vt) \Grad \vt$ is a
continuously differentiable function of the temperature satisfying
\begin{equation} \label{c10}
	0 < \underline{\kappa} \left(1 + \vt^\beta \right) \leq  \kappa(\vt) \leq \Ov{\kappa} \left( 1 + \vt^\beta \right) \ \mbox{for some}\  \beta > 6.
\end{equation}
The restriction $\beta > 6$ is purely technical in the proof of existence of weak solutions. Indeed, if the boundary temperature is constant, one can assume $\beta \geq 3$. Here, the case $\beta = 3$ reflects the effect of radiation.

Finally, we suppose the magnetic diffusivity coefficient $\zeta = \zeta(\vt)$ is a continuously differentiable function of the temperature with
\begin{equation} \label{c11}
	0 < \underline{\zeta}(1 + \vt) \leq \zeta (\vt) \leq \Ov{\zeta}(1 + \vt),\ |\zeta' (\vt) | \leq c
	\ \mbox{for all}\ \vt \geq 0.	
\end{equation}



\section{Relative energy}
\label{r}

A suitable form of the \emph{scaled relative energy} for the compressible MHD system reads
\begin{align}
	E_\eps &\left( \vr, \vt, \vu, \vB \ \Big| r, \Theta, \vU, \vc H \right) =
	\frac{1}{2} \vr |\vu - \vU|^2 + \frac{1}{\eps^2} \frac{1}{2} |\vB - \vc H|^2 \br &+ 	
	\frac{1}{\eps^2} \bigg( \vr e(\vr, \vt) - \Theta \Big( \vr s(\vr, \vt) - r s(r, \Theta) \Big)
	- \Big( e(r, \Theta) - \Theta s(r, \Theta) + \frac{p(r, \Theta)}{r}     \Big)(\vr - r) - r e(r, \Theta) \bigg).
	\label{r6}
\end{align}
In applications, the quantity $(\vr, \vt, \vu, \vB)$ stands for a weak solution of the compressible MHD system while
$(r, \Theta, \vU, \vc H)$ are arbitrary sufficiently smooth functions satisfying the compatibility
conditions
\begin{align}
	r &> 0,\ \Theta > 0 \ \mbox{in}\ [0,T] \times \Ov{\Omega}, \ \Theta|_{\partial \Omega} = \Ov{\vt} + \eps \vtB, \br 
	\vU \cdot \vc n|_{\partial \Omega} &= 0, \ (\S(\Theta, \Grad \vU) \cdot \vc n) \times \vc n|_{\del \Omega} = 0, \br 
	\vc H \times \vc n|_{\del \Omega} &= 0, \ \Div \vc H = 0.
\label{r6a}	
	\end{align}

As a consequence of the hypothesis of thermodynamic stability stated in \eqref{c4}, the energy
\[
E(\vr, \vt, \vu, \vB) = \frac{1}{2} \vr |\vu|^2 + \frac{1}{\eps^2} \vr e(\vr, \vt) + \frac{1}{\eps^2} \frac{1}{2} |\vB|^2
\]
is a strictly convex l.s.c.~function if expressed in the conservative variables $(\vr, S := \vr s(\vr, \vt), \vm := \vr \vu, \vB)$, 
see \cite[Chapter 3, Section 3.1]{FeireislNovotny2022} for details.

\subsection{Relative energy inequality}

The relative energy inequality observed in \cite[Section~3.1.3]{FGKS2023} has a straight-forward extension for the rescaled energy \eqref{r6}, hence this inequality reads in our case
\begin{align}
	& \left[ \intO{ E_\eps \left( \vr, \vt, \vu, \vB \ \Big| r, \Theta, \vU, \vc H \right) } \right]_{t = 0}^{t = \tau}
	\br&	+ \int_0^\tau \intO{ 	\frac{\Theta}{\vt} \left( \mathbb{S}(\vt, \Grad \vu) : \Grad \vu - \frac{1}{\eps^2} \frac{\vc{q}
			(\vt, \Grad \vt) \cdot \Grad \vt}{\vt} + \frac{1}{\eps^2} \zeta(\vt) |\Curl \vB|^2 \right) } \dt \br
	&\quad \leq - \int_0^\tau \intO{ \Big( \vr (\vu - \vU) \otimes (\vu - \vU)   + \frac{1}{\eps^2} p (\vr, \vt) \mathbb{I} - \mathbb{S}(\vt, \Grad \vu) \Big) : \Grad \vU   } \dt \br
	&\quad \quad - \frac{1}{\eps^2} \int_0^\tau \intO{ ( \Curl \vB \times \vB ) \cdot \vU  } \dt
	\br
	&\quad \quad- \int_0^\tau \intO{ \vr \Big(  \partial_t \vU + \vU \cdot \Grad \vU -
		\frac{1}{\eps} \Grad G \Big) \cdot (\vu - \vU)     } \dt \br
	&\quad \quad - \frac{1}{\eps^2} \int_0^\tau \intO{\left( \vr \Big( s(\vr, \vt)  - s(r, \Theta) \Big) \partial_t \Theta + \vr \Big(s (\vr, \vt)  - s(r, \Theta) \Big) \vu \cdot \Grad \Theta + \frac{\vc{q} (\vt, \Grad \vt) }{ \vt} \cdot \Grad \Theta                    \right)      } \dt \br
	 	&\quad \quad + \frac{1}{\eps^2} \int_0^\tau \intO{ \left( \left(1 - \frac{\vr}{r} \right) \partial_t p(r, \Theta) - \frac{\vr}{r} \vu \cdot \Grad p(r,\Theta)      \right)     } \dt \br
	&\quad \quad - \frac{1}{\eps^2} \int_0^\tau \intO{  \Big( \vB \cdot \partial_t \vc H - (\vB \times \vu) \cdot \Curl \vc H - \zeta(\vt) \Curl \vB \cdot \Curl \vc H \Big)    } \dt \br
	&\quad \quad +
	\frac{1}{\eps^2} \int_0^\tau \intO{ \vc H \cdot \partial_t \vc H } \dt \br
		\label{r7}
\end{align}
for any sufficiently regular ``test'' functions $(r, \Theta, \vU, \vc H)$ specified in \eqref{r6a}.


\begin{Remark} \label{Rr1}
	
	By a density argument, the regularity class of functions  $(r, \Theta, \vU, \vc H)$ in \eqref{r6a} can be extended
	to
	\begin{equation} \label{r8}
	r \in W^{1,q}((0,T) \times \Omega), \ (\Theta, \vU, \vc H) \in L^q(0,T; W^{2,q}(\Omega)) \cap W^{1,q}(0,T;
	L^q(\Omega)),
		\end{equation}
where $q$ is large enough such that all terms in \eqref{r7} are well defined.	
	\end{Remark}

Finally, we introduce the definition of weak solutions, existence of which was shown in \cite{FGKS2023}.
\begin{Definition}
    We say that $(\vr, \vt, \vu, \vB)$ is a weak solution of the compressible MHD system \eqref{p1}--\eqref{p4} with boundary conditions \eqref{p7}, \eqref{p72}, and initial data 
    $$\vr(0,\cdot) = \vr_0, \quad \vt(0,\cdot) = \vt_0,\quad \vu(0,\cdot) = \vu_0, \quad \vB(0,\cdot) = \vB_0,$$
if the following holds:\\

For the equation of continuity, we have $\vr \in C_{weak}([0,T] ; L^{\frac53}(\Omega))$, $\vr \geq 0$ a.e.~in $(0,T)\times \Omega$, and the integral identity  
\begin{equation}\label{p1_w_1}
    \int_0^\tau \intO{ \vr \partial_t \varphi+ \vr \vu \cdot \Grad \varphi } \dt =  \left[ \intO{ \vr \varphi } \right]_{t=0}^{t=\tau}
    \end{equation}
holds for all $\varphi \in C^1 ([0,T] \times \Ov\Omega)$ as well as renormalized version of \eqref{p1_w_1} 
\begin{equation}\label{p1_r_1}
    \int_0^\tau \intO{ \xi(\vr) \partial_t \varphi + \xi(\vr) \vu \cdot \Grad \varphi + \left( \xi(\vr) - \xi'(\vr)\vr\right) \Div \vu \, \varphi } \dt =  \left[ \intO{ \xi(\vr) \varphi } \right]_{t=0}^{t=\tau}
    \end{equation}
    holds for all $\varphi \in C^1 ([0,T] \times \Ov\Omega)$ and any $\xi \in C^1(\R)$, $\xi' \in C_c (\R)$.\\
    
For the momentum equation, $\vr \vu \in C_{weak} ([0,T]; L^{\frac{5}{4}}(\Omega; \R^3))$, $\vu \in L^2 (0,T; W^{1,2}(\Omega;\R^3))$, $\vu \cdot \vc{n}|_{\partial \Omega} = 0$, $[\S(\vt, \Grad \vu) \cdot \vc n] \times \vc n|_{\del \Omega} = 0$, and the integral identity  
\begin{equation}\label{p2_w}
    \begin{split}
        &\int_0^\tau \intO{ \vr \vu \cdot \partial_t \vc\varphi + \vr \vu \otimes \vu : \Grad \vc\varphi + p(\vr,\vt) \Div \vc\varphi } \dt
    \\
    &= \int_0^\tau \intO{ \S(\vt, \Grad \vu) : \Grad \vc\varphi } \dt - \int_0^\tau \intO{  \left( \vB \otimes \vB - \frac{1}{2} |\vB|^2 \I \right) : \Grad \vc\varphi } \dt \\ 
    &- \int_0^\tau \intO{ \vr \Grad G \cdot \vc\varphi } \dt + \left[ \intO{ \vr \vu \cdot \vc\varphi } \right]_{t=0}^{t=\tau}
    \end{split}
\end{equation}
holds for any $\vc\varphi \in C([0,T] \times \Ov \Omega ;\ \R^3 )$, $\vc\varphi \cdot \vc{n} |_{\partial\Omega} = 0$.\\

For the induction equation, $\vB \in C_{weak} ([0,T]; L^2(\Omega; \R^3))$, and
    \begin{equation}\label{Bdiv}
    \Div \vB (\tau, \cdot) = 0 \quad \mbox{for any }\tau \in [0,T].        
    \end{equation}
In addition, $$\vB \in L^2 (0,T; H_0(\Omega;\R^3)), \quad H_0 = \{\vc b \in L^2(\Omega): \ \Curl \vc b \in L^2(\Omega), \ \Div \vc b = 0, \ \vc b \times \vc n|_{\del \Omega} = 0\},$$
and the integral identity 
    \begin{equation}
        \begin{split}
            &\int_0^\tau \intO{ \vB \cdot \partial_t \vc\varphi -(\vB \times \vu) \cdot \Curl \vc\varphi - \zeta(\vt) \Curl \vB \cdot \Curl \vc\varphi } \dt\\
            &= \left[ \intO{ \vB \cdot \vc\varphi } \right]_{t=0}^{t=\tau}
        \end{split}
    \end{equation}
holds for all $\vc\varphi \in C^1([0,T) \times \Ov\Omega;\R^3),$ $\Div \vc\varphi =0$, $\vc\varphi \times \vc{n}|_{\partial\Omega} = 0$.\\

For the entropy inequality we have that $\vt \in L^\infty(0,T;L^4(\Omega)),$ $\vt > 0$ a.e.~in $(0,T) \times \Omega$, $\vt - \vt_B \in L^2(0,T; W^{1,2}_0(\Omega)),$ $\log (\vt) \in L^2(0,T; W^{1,2}_0(\Omega)),$ and the integral inequality 
\begin{equation}\label{entropy_weak}
\begin{split}
    &\int_{0}^{\tau}\intO{ \vr s(\vr , \vt ) \partial_t \vc\varphi + \vr s(\vr, \vt) \vu \cdot \Grad \varphi  + \frac{\vc{q}(\vt, \Grad \vt)}{\vt} \cdot \Grad \varphi } \dt \\
    &\leq - \int_{0}^{\tau}\intO{ \frac{1}{\vt} \left(\eps^2 \mathbb{S}(\vt, \Grad \vu) : \Grad \vu - \frac{\vc{q}(\vt, \Grad \vt) \cdot \Grad \vt}{\vt} +  \zeta(\vt) |\Curl \vB |^2 \right) \varphi } \dt
    \\
     &+ \left[ \intO{ \vr s(\vr, \vt) \cdot \vc\varphi } \right]_{t=0}^{t=\tau}
\end{split}
\end{equation}
holds for any $\varphi \in C^1_c ( [0,T] \times \Omega ),$ $\varphi \geq 0.$\\

Moreover, the ballistic energy inequality holds in the sense that 
for any $\psi \in C^1([0,T]\times \Ov\Omega )$, $\psi > 0,$ $\psi|_{\del \Omega} = \vt_B$ the following is satisfied
\begin{equation}
\begin{split}
&\left[ \intO{  \left( \frac{1}{2} \vr |\vu|^2 + \vr e (\vr, \vt)  + \frac{1}{2} | \vB|^2 -  \psi \vr s(\vr,\vt) \right) } \right]_{t=0}^{t=\tau} \\
& + \int_0^\tau \intO{ \frac{\psi}{\vt}\left(\eps^2 \mathbb{S}(\vt, \Grad \vu) : \Grad \vu - \frac{\vc{q}(\vt, \Grad \vt) \cdot \Grad \vt}{\vt} + \zeta(\vt) |\Curl \vB |^2 \right) } \dt
\\
&\leq - \int_0^\tau \intO{ \vr s(\vr,\vt) \partial_t \psi + \vr s(\vr,\vt) \vu  \cdot \Grad \psi + \frac{\vc{q}(\vt, \Grad \vt) }{\vt}  \cdot \Grad \psi } \dt 
\\
&+ \int_0^\tau \intO{ \vr \Grad G \cdot \vu } \dt .
\end{split}
\end{equation}
\end{Definition}



\section{Main result}

Throughout the paper, $C$  stands for a generic constant. For brevity, we sometimes writhe $f\lesssim g $ instead of $f \leq Cg$. 

\begin{Theorem}\label{thm:main}
Let the structural hypotheses be as given above. Let $(\vr_\eps, \vt_\eps, \vu_\eps, \vB_\eps)$ be a solution to the primitive system with boundary conditions
\begin{align*}
\vu_\eps \cdot \vc n = 0, \ [\S(\vt_\eps, \Grad \vu_\eps) \cdot \vc n] \times \vc n = 0, \ \vt_\eps = \Ov{\vt} + \eps \vt_B, \ \Ov{\vt} \in (0,\infty), \ \vB_\eps \times \vc n = 0 \ \text{on} \ \del \Omega.
\end{align*}

Let the initial data be well prepared in the sense that
\begin{align*}
\vr_\eps(0) &= \Ov{\vr} + \eps \vr_{\eps, 0}^1, \ \intO{\vr_{\eps, 0}^1} = 0, \ \vt_\eps(0) = \Ov{\vt} + \eps \vt_{\eps, 0}^1, \ \vu_\eps(0) = \vu_{\eps, 0}, \ \vB_\eps(0) = \vB_{\eps,0}= \Ov{\vB} + \eps \vB_{\eps, 0}^1, \br 
\Ov{\vr} &> 0, \ \Ov{\vt}>0, \ \Ov{\vB} = \Ov{b} \vc e_3 \in \R^3.
\end{align*}
Suppose moreover the compatibility conditions
\begin{align*}
\Curl(\Ov{\vB} \times \vU_0) = 0, 
&& \frac{\del p(\Ov{\vr}, \Ov{\vt})}{\del \vr} \Grad \vr_0^1 + \frac{\del p(\Ov{\vr}, \Ov{\vt})}{\del \vt} \Grad \vt_0^1 = \Ov{\vr} \Grad G - \Curl \vB_0^1 \times \Ov{\vB}.
\end{align*}
If additionally
\begin{align}
\|(\vr_{\eps, 0}^1, \vt_{\eps, 0}^1, \vu_{\eps, 0}, \vB_{\eps, 0}^1)\|_{L^\infty(\Omega)} &\lesssim 1, \br 
\intO{ E_\eps \left( \vr_{\eps, 0}, \vt_{\eps, 0}, \vu_{\eps, 0}, \vB_{\eps, 0} \Big| \Ov{\vr} + \eps \vr_0^1, \Ov{\vt} + \eps \vt_0^1, \vc U_0, \Ov{\vB} + \eps \vB_0^1 \right) } &\to 0, \label{in_en_to_0}
\end{align}
then
\begin{equation}\label{en_to_0}
{\rm ess} \sup_{t \in (0,T)} \intO{ E_\eps \left(\vr_\eps, \vt_\eps, \vu_\eps, \vB_\eps \Big| \Ov{\vr} + \eps \vr^1, \Ov{\vt} + \eps \vt^1, \vc U, \Ov{\vB} + \eps \vB^1 \right) } \to 0,
\end{equation}
where the functions $(\vr^1, \vt^1, \vU, \vB^1)$ are the exact (strong) solution to the limiting modified OBM system \eqref{target} with initial data $(\vr_0^1, \vt_0^1, \vU_0, \vB_0^1)$.
\end{Theorem}

\begin{Remark}
Notice that condition \eqref{in_en_to_0} together with $\|(\vr_{\eps, 0}^1, \vt_{\eps, 0}^1, \vu_{\eps, 0}, \vB_{\eps, 0}^1)\|_{L^\infty(\Omega)} \lesssim 1$ means that 
$\vt^1_{\eps,0},$ $\vr^1_{\eps,0},$ $\vu_{\eps,0}$, $\vB^1_{\eps,0}$ converge strongly in $L^2(\Omega)$ to $\vt^1_{0},$ $\vr^1_{0},$ $\vu_{0}$, $\vB^1_{0}$, respectively. 
\end{Remark}

\begin{Remark}
According to the coercivity properties of the relative energy given below in \eqref{ess_rel}--\eqref{res_rel}, this especially shows
\begin{equation*}
\begin{split}
\frac{\vr_\eps - \Ov{\vr}}{\eps} &\to \vr^1 \quad \text{in} \ L^\infty(0,T;L^1(\Omega)),\\
\frac{\vt_\eps - \Ov{\vt}}{\eps} &\to \vt^1 \quad \text{in} \ L^\infty(0,T;L^1(\Omega)),\\
\sqrt{\vr_\eps} \vu_\eps &\to \sqrt{\Ov{\vr}} \vU \quad \text{in} \ L^\infty(0,T;L^2(\Omega)),\\
\frac{\vB_\eps - \Ov{\vB}}{\eps} &\to \vB^1 \quad \text{in} \ L^\infty(0,T;L^2(\Omega)).
\end{split}
\end{equation*}
\end{Remark}



\section{Basic energy estimates}

We anticipate ill-prepared initial data, meaning
\begin{align*}
\intO{ E_\eps \left( \vr_{\eps, 0}, \vt_{\eps, 0}, \vu_{\eps, 0}, \vB_{\eps, 0} \Big| \Ov{\vr}, \Ov{\vt} + \eps \vt_B, 0, \Ov{\vB} \right) } \lesssim 1.
\end{align*}
This is indeed enough to get uniform bounds on the functions under consideration, and $\vt_B$ is a suitable (we can choose the harmonic) extension of the boundary value $\vt_B$ to the whole of $\Omega$. We use also the first order of the temperature $\Ov{\vt} + \eps \vt_B$ as test function in order to comply with the boundary conditions. For the magnetic field this is not needed since $\vB_\eps \times \vc n = \Ov{\vB} \times \vc n = 0$ on $\del \Omega$.\\

We introduce the essential and residual part of a function $f$ as
\begin{align*}
&M_{\rm ess} = \{(t,x) \in (0,T) \times \Omega: \frac12 \Ov{\vr} \leq \vr_\eps(t,x) \leq 2 \Ov{\vr}, \ \frac12 \Ov{\vt} \leq \vt_\eps(t,x) \leq 2 \Ov{\vt} \},\\
&[f]_{\rm ess} = f\cdot {\bf 1}_{M_{\rm ess}}, \quad [f]_{\rm res} = f - [f]_{\rm ess}.
\end{align*}

The relative energy is in this sense coercive as (see \cite[Section~4.4]{FGKS2023})
\begin{equation}\label{ess_rel}
\left[ E_\eps \left( \vr_\eps, \vt_\eps, \vu_\eps, \vB_\eps \Big| r, \Theta, \vU, \vc H \right) \right]_{\rm ess} \gtrsim \left[ \frac{|\vr_\eps - r|^2}{\eps^2} + \frac{|\vt_\eps - \Theta|^2}{\eps^2} + |\vu_\eps - \vU|^2 + \frac{|\vB_\eps - \vc H|^2}{\eps^2} \right]_{\rm ess},
\end{equation}
\begin{equation}\label{res_rel}
\left[ E_\eps \left( \vr_\eps, \vt_\eps, \vu_\eps, \vB_\eps \Big| r, \Theta, \vU, \vc H \right) \right]_{\rm res} \gtrsim \left[\frac{1}{\eps^2} + \frac{\vr_\eps e(\vr_\eps, \vt_\eps)}{\eps^2} + \frac{\vr_\eps s(\vr_\eps, \vt_\eps)}{\eps^2} + \vr_\eps |\vu_\eps|^2 + \frac{|\vB_\eps|^2}{\eps^2} \right]_{\rm res}.
\end{equation}

Consequently, repeating the steps done in \cite[Section~5.1]{BellaFeireislOschmann2023a} and \cite[Section~2]{Kukuvcka2011}, we find for $(\vr_\eps, \vt_\eps, \vu_\eps, \vB_\eps)$ the uniform estimates
\begin{align}
	\| \vu_\eps \|_{L^2(0,T;W^{1,2} (\Omega)) }  &\lesssim 1, \label{be2} \\
	\frac{1}{\eps^2} \left( \| \Grad \log (\vt_\eps) \|_{L^2((0,T) \times \Omega)}^2 + \| \Grad \vt_\eps^{\frac{\beta}{2}} \|_{L^2((0,T) \times \Omega)}^2 \right) &\lesssim 1.
	\label{be3}
	\end{align}
Additionally,
\begin{equation} \label{be4}
	\frac{1}{\eps^2} \|[1]_{\rm res}]\|_{L^\infty(0,T;L^1(\Omega))} \lesssim 1.
\end{equation}
Moreover,
\begin{align}
	\|\vr_\eps |\vu_\eps|^2\|_{L^\infty(0,T;L^1(\Omega))} &\lesssim 1, \br 
	\left\| \left[ \frac{\vr_\eps - \Ov{\vr}}{\eps} \right]_{\rm ess} \right\|_{L^\infty(0,T; L^2(\Omega))} + \left\| \left[ \frac{\vt_\eps - \Ov{\vt}}{\eps} \right]_{\rm ess} \right\|_{L^\infty(0,T; L^2(\Omega))} &\lesssim 1, \br 
	\frac{1}{\eps^2} \| [\vr_\eps]_{\rm res} \|_{L^\infty(0,T; L^{\frac{5}{3}}(\Omega))}^{\frac{5}{3}} + \frac{1}{\eps^2}  \| [\vt_\eps]_{\rm res} \|_{L^\infty(0,T; L^4(\Omega))}^4	 &\lesssim 1.
\label{be5}	
	\end{align}
Combining \eqref{be3}, \eqref{be4}, and \eqref{be5}, we conclude 
\begin{equation} \label{be6}
	\left\| \frac{\log(\vt_\eps) - \log(\Ov{\vt})}{\eps} \right\|_{L^2(0,T; W^{1,2}(\Omega))} + \left\| \frac{ \vt_\eps - \Ov{\vt} }{\eps} \right\|_{L^2(0,T; W^{1,2}(\Omega))} \lesssim 1.
	\end{equation}
Also, we have the bound on the entropy flux
\begin{equation} \label{be7} 
\left \|\left[ \frac{\kappa (\vt_\eps) }{\vt_\eps} \right]_{\rm res} \frac{\Grad \vt_\eps }{\eps} \right\|_{L^q((0,T) \times \Omega)}  \lesssim 1 \ \mbox{for some}\ q > 1.
\end{equation} 
Finally, we get for the magnetic field the uniform bounds
\begin{align}\label{be8}
\frac{1}{\eps^2} \|\vB_\eps - \Ov{\vB}\|_{L^\infty(0,T; L^2(\Omega))}^2 + \frac{1}{\eps^2} \|\Grad \vB_\eps\|_{L^2((0,T) \times \Omega)}^2 \lesssim 1.
\end{align}







\section{Convergence}
\subsection{Limiting functions}
As a consequence of the uniform bounds obtained, we see that, up to subsequences,
\begin{equation}\label{conv11}
\begin{split}
\vr_\eps &\to \Ov{\vr} \quad \text{in} \ L^\frac{5}{3}(\Omega) \ \text{uniformly in} \ t \in (0,T),\\ 
\vt_\eps &\to \Ov{\vt} \quad \text{in} \ L^2(0,T;W^{1,2}(\Omega)) \cap L^\infty(0,T;L^2(\Omega)), \\ 
\vu_\eps &\weak \vu \quad  \text{weakly in} \ L^2(0,T;W^{1,2}(\Omega)), \\
\vB_\eps &\to \Ov{\vB} \quad \text{in} \ L^2(0,T;W^{1,2}(\Omega)) \cap L^\infty(0,T;L^2(\Omega)).
\end{split}
\end{equation}
Moreover, letting $\eps \to 0$ in the weak formulation of the continuity equation, we obtain
\begin{align*}
\Div \vu = 0, \qquad \vu \cdot \vc n|_{\del \Omega} = 0.
\end{align*}
Next, by the uniform bounds on the essential and residual parts, we see
\providecommand{\RR}{\mathfrak{R}}
\providecommand{\TT}{\mathfrak{T}}
\providecommand{\BB}{\mathfrak{B}}
\begin{equation}\label{w-conv_R}
\left[ \frac{\vr_\eps - \Ov{\vr}}{\eps} \right]_{\rm ess} \weakstar \RR \quad \text{weakly-$\ast$ in} \ L^\infty(0,T;L^2(\Omega)),
\end{equation}
\begin{equation}\label{conv_res}
\left[ \frac{\vr_\eps - \Ov{\vr}}{\eps} \right]_{\rm res} \to 0 \quad \text{in} \ L^\infty(0,T;L^\frac{5}{3}(\Omega)),
\end{equation}
\begin{equation}\label{w-conv_T}
\frac{\vt_\eps - \Ov{\vt}}{\eps} \weak \TT \quad \text{weakly in} \ L^2(0,T;W^{1,2}(\Omega)) \ \text{and weakly-$\ast$ in} \ L^\infty(0,T;L^2(\Omega)),
\end{equation}
\begin{equation}\label{w-conv_B}
\frac{\vB_\eps - \Ov{\vB}}{\eps} \weak \BB \quad \text{weakly in} \ L^2(0,T;W^{1,2}(\Omega)) \ \text{and weakly-$\ast$ in} \ L^\infty(0,T;L^2(\Omega)),
\end{equation}
and in view of the boundary conditions for $\vt_\eps$, we get
\begin{align*}
\TT|_{\del \Omega} = \vt_B.
\end{align*}
Finally, performing the limit in the momentum equation, we have in the sense of distributions
\begin{align*}
\frac{\del p(\Ov{\vr}, \Ov{\vt})}{\del \vr} \Grad \RR + \frac{\del p(\Ov{\vr}, \Ov{\vt})}{\del \vt} \Grad \TT = \Ov{\vr} \Grad G + \Curl \mathfrak{B} \times \Ov{\vB}.
\end{align*}
Again by structural correctness, there exists some function $\mathfrak{A} \in L^2(0,T; W^{1,2}(\Omega))$ with $\intO{\mathfrak{A}} = 0$ such that $\Curl \mathfrak{B} \times \Ov{\vB} = - \Grad \mathfrak{A}$. Hence, it follows additionally that
\begin{align*}
\RR \in L^2(0,T;W^{1,2}(\Omega)).
\end{align*}

\subsection{Relative energy}
To finally show Theorem~\ref{thm:main}, let us use in the relative energy the exact (strong) solution to the limiting modified OB system \eqref{target}.  Namely, we will consider 
\begin{align*}
\intO{ E_\eps \left( \vr_\eps, \vt_\eps, \vu_\eps, \vB_\eps \Big| \Ov{\vr} + \eps \vr^1, \Ov{\vt} + \eps \vt^1, \vc U, \Ov{\vB} + \eps \vB^1 \right) },
\end{align*}
 where $(\vr^1, \vt^1, \vU, \vB^1)$ is the strong solution to \eqref{target}. Keep in mind that 
\begin{align*}
\intO{ E_\eps \left( \vr_{\eps, 0}, \vt_{\eps, 0}, \vu_{\eps, 0}, \vB_{\eps, 0} \Big| \Ov{\vr} + \eps \vr_0^1, \Ov{\vt} + \eps \vt_0^1, \vU_0, \Ov{\vB} + \eps \vB_0^1 \right) } \to 0
\end{align*}
from the assumptions on the initial data given in Theorem~\ref{thm:main}. Let us set 
\begin{align*}
    r_\eps = \Ov{\vr} + \eps \vr^1, \quad \Theta_\eps = \Ov{\vt} + \eps \vt^1, \quad \vc H_\eps = \Ov{\vB} + \eps \vB^1,
\end{align*}
and then we may rewrite the relative energy inequity  \eqref{r7} in the following way
\begin{equation}\label{rel-en-01}
\begin{split}
	& \left[ \intO{ E_\eps \left( \vr_\eps, \vt_\eps, \vu_\eps, \vB_\eps \ \Big| r_\eps, \Theta_\eps, \vU, \vc H_\eps \right) } \right]_{t = 0}^{t = \tau}
	\\ &	+ \int_0^\tau \intO{ 	\frac{\Theta_\eps}{\vt_\eps} \left( \S(\vt_\eps, \Grad \vu_\eps) : \Grad \vu_\eps + \frac{1}{\eps^2} \frac{\kappa(\vt_\eps) |\Grad \vt_\eps|^2}{\vt_\eps} + \frac{1}{\eps^2} \zeta(\vt_\eps) |\Curl \vB_\eps|^2 \right) } \dt \\
	&\quad \leq - \int_0^\tau \intO{ \Big( \vr_\eps (\vu_\eps - \vU) \otimes (\vu_\eps - \vU)   + \frac{1}{\eps^2} p(\vr_\eps, \vt_\eps) \mathbb{I} - \S(\vt_\eps, \Grad \vu_\eps) \Big) : \Grad \vU   } \dt \\
	&\quad \quad - \frac{1}{\eps^2} \int_0^\tau \intO{ ( \Curl \vB_\eps \times \vB_\eps ) \cdot \vU  } \dt \\
	&\quad \quad- \int_0^\tau \intO{ \vr_\eps \Big(  \partial_t \vU + \vU \cdot \Grad \vU -
		\frac{1}{\eps} \Grad G \Big) \cdot (\vu_\eps - \vU) } \dt \\
	&\quad \quad - \frac{1}{\eps^2} \int_0^\tau \intO{\left( \vr_\eps \Big( s(\vr_\eps, \vt_\eps)  - s(r_\eps, \Theta_\eps) \Big) [ \partial_t \Theta_\eps + \vu_\eps \cdot \Grad \Theta_\eps ] - \frac{\kappa(\vt_\eps) \Grad \vt_\eps }{ \vt_\eps} \cdot \Grad \Theta_\eps                    \right) } \dt \\
	 	&\quad \quad + \frac{1}{\eps^2} \int_0^\tau \intO{ \left( \left(1 - \frac{\vr_\eps}{r_\eps} \right) \partial_t p(r_\eps, \Theta_\eps) - \frac{\vr_\eps}{r_\eps} \vu_\eps \cdot \Grad p(r_\eps,\Theta_\eps)      \right)     } \dt
	\\
	&\quad \quad - \frac{1}{\eps^2} \int_0^\tau \intO{  \Big( \vB_\eps \cdot \partial_t \vc H_\eps - (\vB_\eps \times \vu_\eps) \cdot \Curl \vc H_\eps - \zeta(\vt_\eps) \Curl \vB_\eps \cdot \Curl \vc H_\eps \Big)    } \dt \\
	&\quad \quad +
	\frac{1}{\eps^2} \int_0^\tau \intO{ \vc H_\eps \cdot \partial_t \vc H_\eps } \dt .
\end{split}
\end{equation}

Our goal is to rewrite this inequality such that the left-hand side is non-negative, whereas the remainder on the right-hand side will be estimated by the time-integral over the relative energy itself. A Gr\"onwall argument will finally yield the result.

\paragraph{Step 1:} 
Let us rewrite \eqref{rel-en-01} using $\Div \vU=0$, equation \eqref{target}$_2$ (satisfied in a weak sense), and rearranging some terms to obtain
\begin{align}
	& \left[ \intO{ E_\eps \left( \vr_\eps, \vt_\eps, \vu_\eps, \vB_\eps \ \Big| r_\eps, \Theta_\eps, \vU, \vc H_\eps \right) } \right]_{t = 0}^{t = \tau}
	\br&	+ \int_0^\tau \intO{ 	\frac{\Theta_\eps}{\vt_\eps} \left( \S(\vt_\eps, \Grad \vu_\eps) : \Grad \vu_\eps + \frac{1}{\eps^2} \frac{\kappa(\vt_\eps) |\Grad \vt_\eps|^2}{\vt_\eps} + \frac{1}{\eps^2} \zeta(\vt_\eps) |\Curl \vB_\eps|^2 \right) } \dt \br
	&\quad \leq - \frac{1}{\eps} \int_0^\tau \intO{ \vr_\eps \Big( s(\vr_\eps, \vt_\eps)  - s(r_\eps, \Theta_\eps) \Big) [ \partial_t \vt^1 + \vu_\eps \cdot \Grad \vt^1 ] } \dt \br 
	&\quad \quad + \frac{1}{\eps^2} \int_0^\tau \intO{ \frac{\kappa(\vt_\eps) }{ \vt_\eps} \Grad \vt_\eps \cdot \Grad \Theta_\eps } \dt \br 
	&\quad \quad - \int_0^\tau \intO{ [\vr_\eps (\vu_\eps - \vU) \otimes (\vu_\eps - \vU) - \S(\vt_\eps, \Grad \vu_\eps)] : \Grad \vU } \dt \br 
	&\quad \quad - \int_0^\tau \intO{ \frac{\vr_\eps}{\Ov{\vr}} \Big( - \frac{1}{\eps} \Ov{\vr} \Grad G + \vr^1 \Grad G + \Curl \vB^1 \times \vB^1 - \Grad \Pi + \Div \S(\Ov{\vt}, \Grad \vU) \Big) \cdot (\vu_\eps - \vU) } \dt \br
	 	&\quad \quad + \frac{1}{\eps^2} \int_0^\tau \intO{ \left( \left(1 - \frac{\vr_\eps}{r_\eps} \right) \partial_t p(r_\eps, \Theta_\eps) - \frac{\vr_\eps}{r_\eps} \vu_\eps \cdot \Grad p(r_\eps,\Theta_\eps)      \right)     } \dt \br
	&\quad \quad - \frac{1}{\eps^2} \int_0^\tau \intO{ ( \Curl \vB_\eps \times \vB_\eps ) \cdot \vU  } \dt \br
	&\quad \quad - \frac{1}{\eps^2} \int_0^\tau \intO{  \Big( (\vB_\eps - \vc H_\eps) \cdot \partial_t \vc H_\eps - (\vB_\eps \times \vu_\eps) \cdot \Curl \vc H_\eps - \zeta(\vt_\eps) \Curl \vB_\eps \cdot \Curl \vc H_\eps \Big)    } \dt \br 
	&\quad = - \frac{1}{\eps} \int_0^\tau \intO{ \vr_\eps \Big( s(\vr_\eps, \vt_\eps)  - s(r_\eps, \Theta_\eps) \Big) [ \partial_t \vt^1 + \vu_\eps \cdot \Grad \vt^1 ] } \dt \br 
	&\quad \quad + \frac{1}{\eps^2} \int_0^\tau \intO{ \frac{\kappa(\vt_\eps) }{ \vt_\eps} \Grad \vt_\eps \cdot \Grad \Theta_\eps } \dt \br 
	&\quad \quad - \int_0^\tau \intO{ [\vr_\eps (\vu_\eps - \vU) \otimes (\vu_\eps - \vU) - \S(\vt_\eps, \Grad \vu_\eps)] : \Grad \vU } \dt \br 
	&\quad \quad - \int_0^\tau \intO{ \frac{\vr_\eps}{\Ov{\vr}} \Big( - \frac{1}{\eps} \Ov{\vr} \Grad G + \vr^1 \Grad G - \Grad \Pi \Big) \cdot (\vu_\eps - \vU) } \dt \br 
	&\quad \quad + \int_0^\tau \intO{ \Big( \frac{\vr_\eps}{\Ov{\vr}} - 1 \Big) \Big( \Div \S(\Ov{\vt}, \Grad \vU) + \Curl \vB^1 \times \vB^1 \Big) \cdot (\vU - \vu_\eps) } \dt \br 
	&\quad \quad + \int_0^\tau \intO{ \S(\Ov{\vt}, \Grad \vU) : \Grad (\vu_\eps - \vU) } \dt \br 
	&\quad \quad + \frac{1}{\eps^2} \int_0^\tau \intO{ \left( \left(1 - \frac{\vr_\eps}{r_\eps} \right) \partial_t p(r_\eps, \Theta_\eps) - \frac{\vr_\eps}{r_\eps} \vu_\eps \cdot \Grad p(r_\eps,\Theta_\eps)      \right)     } \dt \br 
    &\quad \quad - \frac{1}{\eps} \int_0^\tau \intO{ ( \Curl \vB^1 \times \Ov{\vB} ) \cdot (\vU - \vu_\eps) } \dt \br
	&\quad \quad + \frac{1}{\eps^2} \int_0^\tau \intO{ ( \Curl \vc H_\eps \times \vc H_\eps ) \cdot (\vU - \vu_\eps) } \dt \br
	&\quad \quad - \frac{1}{\eps^2} \int_0^\tau \intO{ ( \Curl \vB_\eps \times \vB_\eps ) \cdot \vU  } \dt \br
	&\quad \quad - \frac{1}{\eps^2} \int_0^\tau \intO{  \Big( (\vB_\eps - \vc H_\eps) \cdot \partial_t \vc H_\eps - (\vB_\eps \times \vu_\eps) \cdot \Curl \vc H_\eps - \zeta(\vt_\eps) \Curl \vB_\eps \cdot \Curl \vc H_\eps \Big)    } \dt . \label{first}
\end{align}

We see that we can get immediately rid of some terms: by strong convergence $\vr_\eps \to \Ov{\vr}$ in $L^\infty(0,T; L^\frac{5}{3}(\Omega))$ and weak convergence $\vu_\eps \weak \vu$ in $L^2(0,T; L^6(\Omega))$, together with $\Div \vU = 0 = \Div \vu$, we have
\begin{align*}
    &\int_0^\tau \intO{ \frac{\vr_\eps}{\Ov{\vr}} \Grad \Pi \cdot (\vu_\eps - \vU) } \dt \br 
    &\quad + \int_0^\tau \intO{ \Big( \frac{\vr_\eps}{\Ov{\vr}} - 1 \Big) \Big( \Div \S(\Ov{\vt}, \Grad \vU) + \Curl \vB^1 \times \vB^1 \Big) \cdot (\vU - \vu_\eps) } \dt \br 
    &= \int_0^\tau \intO{ \Grad \Pi \cdot (\vu - \vU) } \dt + O(\eps) = O(\eps),
\end{align*}
where $O(\eps)$ is a generic error that vanishes as $\eps \to 0$.

\paragraph{Step 2:} 
Let us focus on the very last term of \eqref{first}. 
As $\Curl(\Ov{\vB} \times \vU) = 0$, we see that $\vc H_\eps$ satisfies the induction equation. As also $(\vB_\eps - \vc H_\eps) \times \vc n |_{\del \Omega} = 0$, integration by parts gives
\begin{align*}
& - \frac{1}{\eps^2} \int_0^\tau \intO{  \Big( (\vB_\eps - \vc H_\eps) \cdot \partial_t \vc H_\eps - (\vB_\eps \times \vu_\eps) \cdot \Curl \vc H_\eps - \zeta(\vt_\eps) \Curl \vB_\eps \cdot \Curl \vc H_\eps \Big)    } \dt \br 
&= \frac{1}{\eps^2} \int_0^\tau \intO{  \Big( (\vB_\eps - \vc H_\eps) \cdot \Curl (\vc H_\eps \times \vc U) + (\vB_\eps \times \vu_\eps) \cdot \Curl \vc H_\eps } \dt \br 
 &\quad + \frac{1}{\eps^2} \int_0^\tau \intO{ (\vB_\eps - \vc H_\eps) \cdot \Curl (\zeta(\Ov{\vt}) \Curl \vc H_\eps) + \zeta(\vt_\eps) \Curl \vB_\eps \cdot \Curl \vc H_\eps \Big)    } \dt \br 
 &= \frac{1}{\eps^2} \int_0^\tau \intO{  \Big( (\vB_\eps - \vc H_\eps) \cdot \Curl (\vc H_\eps \times \vc U) + (\vB_\eps \times \vu_\eps) \cdot \Curl \vc H_\eps } \dt \br 
 &\quad + \frac{1}{\eps^2} \int_0^\tau \intO{ \zeta(\Ov{\vt}) \Curl (\vB_\eps - \vc H_\eps) \cdot \Curl \vc H_\eps + \zeta(\vt_\eps) \Curl \vB_\eps \cdot \Curl \vc H_\eps \Big)    } \dt .
\end{align*}

Gathering terms not containing magnetic dissipation and using \cite[Equation (4.6)]{FGKS2023}, we get
\begin{align*}
&\frac{1}{\eps^2} \intO{ \Curl \vc H_\eps \times \vc H_\eps \cdot  (\vU - \vu_\eps) - \Curl \vB_\eps \times \vB_\eps \cdot \vU } \br 
&\quad + \frac{1}{\eps^2} \intO{ (\vB_\eps - \vc H_\eps) \cdot \Curl (\vc H_\eps \times \vU) + (\vB_\eps \times \vu_\eps) \cdot \Curl \vc H_\eps} \br 
&= \frac{1}{\eps^2} \intO{ (\vU \times (\vB_\eps - \vc H_\eps)) \cdot \Curl(\vB_\eps - \vc H_\eps) } + \frac{1}{\eps^2} \intO{ ((\vU - \vu_\eps) \times (\vB_\eps - \vc H_\eps)) \cdot \Curl \vc H_\eps }.
\end{align*}
Hence, we  find that 
\begin{align}\label{rel-en-03}
	& \left[ \intO{ E_\eps\left( \vr_\eps, \vt_\eps, \vu_\eps, \vB_\eps \ \Big| r_\eps, \Theta_\eps, \vU, \vc H_\eps \right) } \right]_{t = 0}^{t = \tau}
	\br&	+ \int_0^\tau \intO{ 	\frac{\Theta_\eps}{\vt_\eps} \left( \S(\vt_\eps, \Grad \vu_\eps) : \Grad \vu_\eps + \frac{1}{\eps^2} \frac{\kappa(\vt_\eps) |\Grad \vt_\eps|^2}{\vt_\eps} + \frac{1}{\eps^2} \zeta(\vt_\eps) |\Curl \vB_\eps|^2 \right) } \dt \br
	&\quad \leq - \frac{1}{\eps} \int_0^\tau \intO{ \vr_\eps \Big( s(\vr_\eps, \vt_\eps)  - s(r_\eps, \Theta_\eps) \Big) [ \partial_t \vt^1 + \vu_\eps \cdot \Grad \vt^1 ] } \dt \br 
	&\quad \quad + \frac{1}{\eps^2} \int_0^\tau \intO{ \frac{\kappa(\vt_\eps) }{ \vt_\eps} \Grad \vt_\eps \cdot \Grad \Theta_\eps } \dt \br 
	&\quad \quad - \int_0^\tau \intO{ [\vr_\eps (\vu_\eps - \vU) \otimes (\vu_\eps - \vU) ] : \Grad \vU } \dt \br 
	&\quad \quad + \int_0^\tau \intO{ \frac{\vr_\eps}{\Ov{\vr}} \Big( \vr^1 \Grad G - \frac{1}{\eps} \Ov{\vr} \Grad G \Big) \cdot (\vU - \vu_\eps) } \dt \br 
	&\quad \quad + \int_0^\tau \intO{ \S(\Ov{\vt}, \Grad \vU) : \Grad (\vu_\eps - \vU) } \dt + \int_0^\tau \intO{ \S(\vt_\eps, \Grad \vu_\eps) : \Grad \vU } \dt \br 
	 	&\quad \quad + \frac{1}{\eps^2} \int_0^\tau \intO{ \left( \left(1 - \frac{\vr_\eps}{r_\eps} \right) \partial_t p(r_\eps, \Theta_\eps) - \frac{\vr_\eps}{r_\eps} \vu_\eps \cdot \Grad p(r_\eps,\Theta_\eps)      \right)     } \dt \br 
        &\quad \quad - \frac{1}{\eps} \int_0^\tau \intO{ ( \Curl \vB^1 \times \Ov{\vB} ) \cdot (\vU - \vu_\eps) } \dt \br
	&\qquad + \frac{1}{\eps^2} \intO{ (\vU \times (\vB_\eps - \vc H_\eps)) \cdot \Curl(\vB_\eps - \vc H_\eps) } \br 
    &\qquad + \frac{1}{\eps^2} \intO{ ((\vU - \vu_\eps) \times (\vB_\eps - \vc H_\eps)) \cdot \Curl \vc H_\eps } \dt \br 
 &\qquad + \frac{1}{\eps^2} \int_0^\tau \intO{ \zeta(\Ov{\vt}) \Curl (\vB_\eps - \vc H_\eps) \cdot \Curl \vc H_\eps + \zeta(\vt_\eps) \Curl \vB_\eps \cdot \Curl \vc H_\eps     } \dt + O(\eps) .
\end{align}
Lastly, inserting the rewritten heat equation \eqref{target}$_3$ satisfied by $\vt^1$, where $\alpha = \alpha(\Ov{\vr}, \Ov{\vt})$, $c_p = c_p(\Ov{\vr}, \Ov{\vt})$, specifically,
\begin{align*}
\del_t \vt^1 + \vU \cdot \Grad \vt^1 = \frac{\kappa(\Ov{\vt})}{\Ov{\vr} c_p} \Delta \vt^1 + \frac{\Ov{\vt} \alpha}{c_p} \vU \cdot \Grad G + \frac{\Ov{\vt} \alpha}{\Ov{\vr} c_p} \del_t \chi - \frac{\Ov{\vt} \alpha}{\Ov{\vr} c_p} (\del_t A + \vU \cdot \Grad A),
\end{align*}
we may rewrite the first term on the RHS of \eqref{rel-en-03} to get
\begin{align}\label{ei1}
	& - \frac{1}{\eps} \int_0^\tau \intO{ \vr_\eps \Big( s(\vr_\eps, \vt_\eps)  - s(r_\eps, \Theta_\eps) \Big) [ \partial_t \vt^1 + \vu_\eps \cdot \Grad \vt^1 ] } \dt \br
	&= - \frac{1}{\eps} \int_0^\tau \intO{ \vr_\eps \Big( s(\vr_\eps, \vt_\eps)  - s(r_\eps, \Theta_\eps) \Big) \left[ \frac{\kappa(\Ov{\vt})}{\Ov{\vr} c_p} \Delta \vt^1 + \frac{\Ov{\vt} \alpha}{c_p} \vU \cdot \Grad G + \frac{\Ov{\vt} \alpha}{\Ov{\vr} c_p} \del_t \chi \right] } \dt \br 
    &\qquad - \frac{1}{\eps} \int_0^\tau \intO{ \vr_\eps \Big( s(\vr_\eps, \vt_\eps)  - s(r_\eps, \Theta_\eps) \Big) [ (\vu_\eps - \vU) \cdot \Grad \vt^1 ] } \dt \br 
    &\qquad + \frac{1}{\eps} \int_0^\tau \intO{ \vr_\eps \Big( s(\vr_\eps, \vt_\eps)  - s(r_\eps, \Theta_\eps) \Big) \left[ \frac{\Ov{\vt} \alpha}{\Ov{\vr} c_p} (\del_t A + \vU \cdot \Grad A) \right] } \dt .
\end{align}

\paragraph{Step 3:} As for $s(\vr_\eps, \vt_\eps)$ in \eqref{ei1}, we write
\footnote{One needs to add references(\cite{FeireislNovotny2009singlim}) and write it in a more proper form, integral. Answer: better? }
\begin{align*}
&\int_0^\tau \intO{ \psi \vr_\eps \frac{s(\vr_\eps, \vt_\eps) - s(r_\eps, \Theta_\eps)}{\eps} } \dt = \int_0^\tau \intO{ \psi \vr_\eps \frac{s(\vr_\eps, \vt_\eps) - s(\Ov{\vr}, \Ov{\vt})}{\eps} + \psi \vr_\eps \frac{s(\Ov{\vr}, \Ov{\vt}) - s(r_\eps, \Theta_\eps)}{\eps} } \dt \\
&= \int_0^\tau \intO{ \psi \vr_\eps \left( \del_\vr s(\Ov{\vr}, \Ov{\vt}) \frac{\vr_\eps - \Ov{\vr}}{\eps} + \del_\vt s(\Ov{\vr}, \Ov{\vt}) \frac{\vt_\eps - \Ov{\vt}}{\eps} - \del_\vr s(\Ov{\vr}, \Ov{\vt}) \vr^1 - \del_\vt s(\Ov{\vr}, \Ov{\vt}) \vt^1 \right) } \dt + O(\eps) \\
&= \int_0^\tau \intO{ \psi \Ov{\vr} \left( \del_\vr s(\Ov{\vr}, \Ov{\vt}) (\mathfrak{R} - \vr^1) + \del_\vt s(\Ov{\vr}, \Ov{\vt}) (\mathfrak{T} - \vt^1) \right) } \dt + O(\eps),
\end{align*}
where $\psi$ represents each of the factors in $[...]$ in \eqref{ei1}.

For the term containing $p(r_\eps, \Theta_\eps)$ in \eqref{rel-en-03}, we use the Boussinesq relation \eqref{OB} pointwise to write
\begin{align*}
\del_t p(r_\eps, \Theta_\eps) &= \eps( \del_\vr p(r_\eps, \Theta_\eps) \del_t \vr^1 + \del_\vt p(r_\eps, \Theta_\eps) \del_t \vt^1 ) \br 
&= \eps [ (\del_\vr p(r_\eps, \Theta_\eps) - \del_\vr p(\Ov{\vr}, \Ov{\vt})) \del_t \vr^1 + (\del_\vt p(r_\eps, \Theta_\eps) - \del_\vt p(\Ov{\vr}, \Ov{\vt})) \del_t \vt^1 ] \\
&\quad + \eps [ \del_\vr p(\Ov{\vr}, \Ov{\vt}) \del_t \vr^1 + \del_\vt p(\Ov{\vr}, \Ov{\vt} ) \del_t \vt^1] \\
&= \eps^2 [ \del_\vr^2 p(\Ov{\vr}, \Ov{\vt}) \vr^1 \del_t \vr^1 + \del_{\vr \vt}^2 p(\Ov{\vr}, \Ov{\vt}) \del_t (\vr^1 \vt^1) + \del_\vt^2 p(\Ov{\vr}, \Ov{\vt}) \vt^1 \del_t \vt^1 ] + \eps \del_t (\chi - A) + O(\eps),
\end{align*}
and similarly,
\begin{align*}
\Grad p(r_\eps, \Theta_\eps) = \eps^2 [ \del_\vr^2 p(\Ov{\vr}, \Ov{\vt}) \vr^1 \Grad \vr^1 + \del_{\vr \vt}^2 p(\Ov{\vr}, \Ov{\vt}) \Grad (\vr^1 \vt^1) + \del_{\vt}^2 p(\Ov{\vr}, \Ov{\vt}) \vt^1 \Grad \vt ] + \eps \Grad (\Ov{\vr} G - A) + O(\eps).
\end{align*}

Thus, we find
\begin{align*}
&\frac{1}{\eps^2} \int_0^\tau \intO{ \Big( 1 - \frac{\vr_\eps}{r_\eps} \Big) \del_t p(r_\eps, \Theta_\eps) - \frac{\vr_\eps}{r_\eps} \vu_\eps \cdot \Grad p(r_\eps, \Theta_\eps) } \dt \\
&= \int_0^\tau \intO{ \Big( 1 - \frac{\vr_\eps}{r_\eps} \Big) ( \del_\vr^2 p(\Ov{\vr}, \Ov{\vt}) \vr^1 \del_t \vr^1 + \del_{\vr \vt}^2 p(\Ov{\vr}, \Ov{\vt}) \del_t (\vr^1 \vt^1) + \del_\vt^2 p(\Ov{\vr}, \Ov{\vt}) \vt^1 \del_t \vt^1 + \frac{1}{\eps} (\del_t \chi - \del_t A) ) } \dt \br 
&\quad - \int_0^\tau \intO{ \frac{\vr_\eps}{r_\eps} \vu_\eps \cdot (\del_\vr^2 p(\Ov{\vr}, \Ov{\vt}) \vr^1 \Grad \vr^1 + \del_{\vr \vt}^2 p(\Ov{\vr}, \Ov{\vt}) \Grad (\vr^1 \vt^1) + \del_{\vt}^2 p(\Ov{\vr}, \Ov{\vt}) \vt^1 \Grad \vt ) } \dd t \br 
&\quad - \int_0^\tau \intO{ \frac{1}{\eps} \frac{\vr_\eps}{r_\eps} \vu_\eps \cdot ( \Grad \Ov{\vr} G - \Grad A) } \dt + O(\eps) \br 
&= - \frac{1}{\eps} \int_0^\tau \intO{ \frac{\vr_\eps}{r_\eps} \vu_\eps \cdot (\Grad \Ov{\vr} G - \Grad A) } \dd t - \int_0^\tau \intO{ \frac{1}{\Ov{\vr}} (\vr^1 - \mathfrak{R}) \del_t A } \dt + O(\eps),
\end{align*}
the last line coming from strong convergence $\vr_\eps/r_\eps \to 1$, say in $L^\infty(0,T;L^{\frac{5}{3}}(\Omega))$, as $\eps \to 0$, weak convergence $\vu_\eps \weak \vu$ in $L^2(0,T;L^6(\Omega))$, as well as $\Div \vu = 0$, and the fact that for $\del_t (\chi - A)$, we have
\begin{align*}
&\frac{1}{\eps} \int_0^\tau \intO{ \Big( 1 - \frac{\vr_\eps}{r_\eps} \Big) \del_t (\chi - A) } \dt = \int_0^\tau \intO{ \frac{1}{r_\eps} \frac{(\Ov{\vr} + \eps \vr^1) - \vr_\eps}{\eps} \del_t (\chi - A) } \dt \\
&= \int_0^\tau \intO{ \frac{1}{\Ov{\vr}} (\vr^1 - \mathfrak{R}) \del_t (\chi - A) } \dt + O(\eps) = - \int_0^\tau \intO{ \frac{1}{\Ov{\vr}} (\vr^1 - \mathfrak{R}) \del_t A } \dt + O(\eps)
\end{align*}
due to $\int_\Omega \vr^1 \dx = \int_\Omega \mathfrak{R} \dx = 0$ and $\chi$ is spatially homogeneous. Thus, putting also diffusive terms to the left hand site of \eqref{rel-en-03},

\begin{align}\label{ei2}
	& \left[ \intO{ E_\eps\left( \vr_\eps, \vt_\eps, \vu_\eps, \vB_\eps \ \Big| r_\eps, \Theta_\eps, \vU, \vc H_\eps \right) } \right]_{t = 0}^{t = \tau}
	\br&	+ \int_0^\tau \intO{ 	\frac{\Theta_\eps}{\vt_\eps} \left( \S(\vt_\eps, \Grad \vu_\eps) : \Grad \vu_\eps + \frac{1}{\eps^2} \frac{\kappa(\vt_\eps) |\Grad \vt_\eps|^2}{\vt_\eps} + \frac{1}{\eps^2} \zeta(\vt_\eps) |\Curl \vB_\eps|^2 \right) } \dt \br 
	& - \int_0^\tau \intO{ \S(\Ov{\vt}, \Grad \vU) : \Grad (\vu_\eps - \vU) } \dt - \int_0^\tau \intO{ \S(\vt_\eps, \Grad \vu_\eps) : \Grad \vU } \dt \br 
	&- \frac{1}{\eps^2} \int_0^\tau \intO{ \frac{\kappa(\vt_\eps) }{ \vt_\eps} \Grad \vt_\eps \cdot \Grad \Theta_\eps } \dt \br 
	&- \frac{1}{\eps^2} \int_0^\tau \intO{ \zeta(\Ov{\vt}) \Curl (\vB_\eps - \vc H_\eps) \cdot \Curl \vc H_\eps + \zeta(\vt_\eps) \Curl \vB_\eps \cdot \Curl \vc H_\eps     } \dt \br 
	&\quad \leq R, 
\end{align}
where the remainder is given by
\begin{align*}
    R &= - \int_0^\tau \intO{ \Ov{\vr} \Big( \del_\vr s(\Ov{\vr}, \Ov{\vt}) (\mathfrak{R} - \vr^1) + \del_\vt s(\Ov{\vr}, \Ov{\vt}) (\mathfrak{T} - \vt^1) \Big) \left[ \frac{\kappa(\Ov{\vt})}{\Ov{\vr} c_p} \Delta \vt^1 + \frac{\Ov{\vt} \alpha}{c_p} \vU \cdot \Grad G + \frac{\Ov{\vt} \alpha}{\Ov{\vr} c_p} \del_t \chi \right] } \dt \br 
    &\qquad - \frac{1}{\eps} \int_0^\tau \intO{ \vr_\eps \Big( s(\vr_\eps, \vt_\eps)  - s(r_\eps, \Theta_\eps) \Big) [ (\vu_\eps - \vU) \cdot \Grad \vt^1 ] } \dt \br 
    &\qquad + \int_0^\tau \intO{ \Ov{\vr} \Big( \del_\vr s(\Ov{\vr}, \Ov{\vt}) (\mathfrak{R} - \vr^1) + \del_\vt s(\Ov{\vr}, \Ov{\vt}) (\mathfrak{T} - \vt^1) \Big) \left[ \frac{\Ov{\vt} \alpha}{\Ov{\vr} c_p} (\del_t A + \vU \cdot \Grad A) \right] } \br 
	&\quad \quad - \int_0^\tau \intO{ [\vr_\eps (\vu_\eps - \vU) \otimes (\vu_\eps - \vU) ] : \Grad \vU } \dt \br 
	&\qquad + \int_0^\tau \intO{ \frac{\vr_\eps}{\Ov{\vr}} \Big( \vr^1 \Grad G - \frac{1}{\eps} \Ov{\vr} \Grad G \Big) \cdot (\vU - \vu_\eps) } \dt \br 
	 	&\quad \quad - \frac{1}{\eps} \int_0^\tau \intO{ \frac{\vr_\eps}{r_\eps} \vu_\eps \cdot (\Ov{\vr} \Grad G - \Grad A)  } \dt \br
        &\quad \quad - \int_0^\tau \intO{ \frac{1}{\Ov{\vr}} (\vr^1 - \mathfrak{R}) \del_t A  } \dt \br
        &\quad \quad - \frac{1}{\eps} \int_0^\tau \intO{ ( \Curl \vB^1 \times \Ov{\vB} ) \cdot (\vU - \vu_\eps) } \dt \br
	&\qquad + \frac{1}{\eps^2} \int_0^\tau \intO{ (\vU \times (\vB_\eps - \vc H_\eps)) \cdot \Curl(\vB_\eps - \vc H_\eps) } \dd t \br 
	&\qquad + \frac{1}{\eps^2} \int_0^\tau \intO{ ((\vU - \vu_\eps) \times (\vB_\eps - \vc H_\eps)) \cdot \Curl \vc H_\eps } \dd t + O(\eps).
\end{align*}

\paragraph{Step 4:}
By the form of relative energy we get that 
\begin{equation*}
\frac{1}{\eps} \int_0^\tau \intO{ \vr_\eps (s(\vr_\eps, \vt_\eps) - s(r_\eps, \Theta_\eps)) (\vu_\eps - \vU) \cdot \Grad \vt^1 } \dd t \lesssim \int_0^\tau E_\eps\left( \vr_\eps, \vt_\eps, \vu_\eps, \vB_\eps \ \Big| r_\eps, \Theta_\eps, \vU, \vc H_\eps \right) \dt,
\end{equation*}
\begin{equation*}
-\int_0^\tau \intO{ \vr_\eps (\vu_\eps - \vU) \otimes (\vu_\eps - \vU) : \Grad \vU } \dd t \lesssim \int_0^\tau E_\eps\left( \vr_\eps, \vt_\eps, \vu_\eps, \vB_\eps \ \Big| r_\eps, \Theta_\eps, \vU, \vc H_\eps \right) \dd t,
\end{equation*}
where the very first inequality can be proven similarly to \cite[Section~4.2.4, eqs. (4.46)--(4.48)]{FeireislNovotny2022}.\\

The same way as in \cite{BellaFeireislOschmann2023a}, we find also
\begin{align*}
&-\frac{1}{\eps} \int_0^\tau \intO{ \frac{\vr_\eps}{r_\eps} \vu_\eps \cdot (\Ov{\vr} \Grad G - \Grad A) } \dd t \\
&= - \frac{1}{\eps} \int_0^\tau \intO{ \vr_\eps \vu_\eps \cdot \Big( \Grad G - \frac{1}{\Ov{\vr}} \Grad A \Big) } \dt + \frac{1}{\eps} \int_0^\tau \intO{ \Big( 1 - \frac{\Ov{\vr}}{r_\eps} \Big) \vr_\eps \vu_\eps \cdot \Big( \Grad G - \frac{1}{\Ov{\vr}} \Grad A \Big) } \dt \\
&= - \frac{1}{\eps} \int_0^\tau \intO{ \vr_\eps \vu_\eps \cdot \Big( \Grad G - \frac{1}{\Ov{\vr}} \Grad A \Big) \dd x \dd t + \int_0^\tau \int_\Omega \vr^1 \vu \cdot \Big( \Grad G - \frac{1}{\Ov{\vr}} \Grad A \Big) } \dd t + O(\eps).
\end{align*}
Additionally, since \eqref{w-conv_R}--\eqref{w-conv_B} hold, we find from the weak form of the momentum equation that
\begin{align*}
\del_\vr p(\Ov{\vr}, \Ov{\vt}) \Grad \mathfrak{R} + \del_\vt p(\Ov{\vr}, \Ov{\vt}) \Grad \mathfrak{T} = \Ov{\vr} \Grad G - \Grad \mathfrak{A} 
\end{align*}
in the sense of distributions. That combined with \eqref{OB} leads to
\begin{align}\label{BoussDiff}
\del_\vr p(\Ov{\vr}, \Ov{\vt}) \Grad (\mathfrak{R} - \vr^1) + \del_\vt p(\Ov{\vr}, \Ov{\vt}) \Grad (\mathfrak{T} - \vt^1) + \Grad(\mathfrak{A} - A) = 0
\end{align}
in the sense of distributions. Hence, using that by definition $- \Grad A = \Curl \vB^1 \times \Ov{\vB}$ and that $\Div \vU = 0$, the remainder becomes
\begin{align*}
R &= - \int_0^\tau \intO{ \Ov{\vr} \Big( \del_\vr s(\Ov{\vr}, \Ov{\vt}) (\mathfrak{R} - \vr^1) + \del_\vt s(\Ov{\vr}, \Ov{\vt}) (\mathfrak{T} - \vt^1) \Big) \left[ \frac{\kappa(\Ov{\vt})}{\Ov{\vr} c_p} \Delta \vt^1 + \frac{\Ov{\vt} \alpha}{c_p} \vU \cdot \Grad G + \frac{\Ov{\vt} \alpha}{\Ov{\vr} c_p} \del_t \chi \right] } \dt \br 
    &\qquad + \int_0^\tau \intO{ \Ov{\vr} \Big( \del_\vr s(\Ov{\vr}, \Ov{\vt}) (\mathfrak{R} - \vr^1) + \del_\vt s(\Ov{\vr}, \Ov{\vt}) (\mathfrak{T} - \vt^1) \Big) \left[ \frac{\Ov{\vt} \alpha}{\Ov{\vr} c_p} (\del_t A + \vU \cdot \Grad A) \right] } \br 
	&\qquad + \int_0^\tau \intO{ \frac{\vr_\eps}{\Ov{\vr}} \vr^1 \Grad G \cdot (\vU - \vu_\eps) } \dt - \frac{1}{\eps} \int_0^\tau \intO{ \vr_\eps \Grad G \cdot \vU  } \dt \br
    &\qquad + \int_0^\tau \intO{ \vr^1 \vu \cdot \Grad G } \dd t - \int_0^\tau \intO{ \frac{1}{\Ov{\vr}} (\vr^1 - \mathfrak{R}) \del_t A  } \dt \br
        &\quad \quad + \frac{1}{\eps} \int_0^\tau \intO{ \Big( \frac{\vr_\eps}{\Ov{\vr}} - 1 \Big) \vu_\eps \cdot \Grad A \dd x \dd t - \int_0^\tau \int_\Omega \vr^1 \vu \cdot \frac{1}{\Ov{\vr}} \Grad A } \dd t \br
	&\qquad + \frac{1}{\eps^2} \int_0^\tau \intO{ (\vU \times (\vB_\eps - \vc H_\eps)) \cdot \Curl(\vB_\eps - \vc H_\eps) } \dd t \br 
	&\qquad + \frac{1}{\eps^2} \int_0^\tau \intO{ ((\vU - \vu_\eps) \times (\vB_\eps - \vc H_\eps)) \cdot \Curl \vc H_\eps } \dd t \\
    &\qquad + \int_0^\tau \intO{ E_\eps\left( \vr_\eps, \vt_\eps, \vu_\eps, \vB_\eps \ \Big| r_\eps, \Theta_\eps, \vU, \vc H_\eps \right) } \dt + O(\eps). 
\end{align*}
By strong convergence $\vr_\eps \to \Ov{\vr}$ in $L^\infty(0,T;L^{\frac{5}{3}}(\Omega))$, we may replace one by the other in the third integral of the above. Using moreover $\Div \vU = 0$, we may write
\begin{align*}
- \frac{1}{\eps} \int_0^\tau \intO{ \vr_\eps \Grad G \cdot \vU } \dt = \int_0^\tau \intO{ \frac{\Ov{\vr} - \vr_\eps}{\eps} \Grad G \cdot \vU } \dt = - \int_0^\tau \intO{ \mathfrak{R} \Grad G \cdot \vU } \dt + O(\eps).
\end{align*}
By the above, we rewrite using $\Div \vU = 0$, weak convergence $\vu_\eps \weak \vu$ in $L^2(0,T; L^6(\Omega))$, and strong convergence $\vr_\eps \to \Ov{\vr}$ in $L^\infty(0,T; L^\frac{5}{3}(\Omega))$
\begin{align}
    &\int_0^\tau \intO{ \frac{\vr_\eps}{\Ov{\vr}} \vr^1 \Grad G \cdot (\vU - \vu_\eps) } \dt - \int_0^\tau \intO{ \frac{1}{\eps} \vr_\eps \Grad G \cdot \vU } \dt + \int_0^\tau \intO{ \vr^1 \vu \cdot \Grad G } \dt \br 
    &= \int_0^\tau \intO{ \Big( \frac{\vr_\eps}{\Ov{\vr}} - 1 \Big) \vr^1 \Grad G \cdot (\vU - \vu_\eps) } \dt + \int_0^\tau \intO{ \Big( \vr^1 - \frac{\vr_\eps - \Ov{\vr}}{\eps} \Big) \Grad G \cdot \vU } \dt \br 
    &\qquad + \int_0^\tau \intO{ \vr^1 \Grad G \cdot (\vu - \vu_\eps) } \dt \br 
    &= \int_0^\tau \intO{ ( \vr^1 - \mathfrak{R} ) \Grad G \cdot \vU } \dt + O(\eps). 
\end{align}

Similarly for $A$, by $\Div \vU = 0$,
\begin{align}
    &\frac{1}{\eps} \int_0^\tau \intO{ \Big( \frac{\vr_\eps}{\Ov{\vr}} - 1 \Big) \vu_\eps \cdot \Grad A } \dt - \int_0^\tau \intO{ \vr^1 \vu \cdot\frac{1}{\Ov{\vr}} \Grad A } \dt \br 
    &=\frac{1}{\eps} \int_0^\tau \intO{ \Big( \frac{\vr_\eps}{\Ov{\vr}} - 1 \Big) (\vu_\eps - \vU) \cdot \Grad A } \dt + \int_0^\tau \intO{ \frac{1}{\Ov{\vr}} \frac{\vr_\eps - \Ov{\vr}}{\eps} \vU \cdot \Grad A } \dt - \int_0^\tau \intO{ \vr^1 \vu \cdot\frac{1}{\Ov{\vr}} \Grad A } \dt \br 
    &=\frac{1}{\eps} \int_0^\tau \intO{ \Big( \frac{\vr_\eps}{\Ov{\vr}} - 1 \Big) (\vu_\eps - \vU) \cdot \Grad A } \dt + \int_0^\tau \intO{ \frac{1}{\Ov{\vr}} \Big( \frac{\vr_\eps - \Ov{\vr}}{\eps} - \vr^1 \Big) \vU \cdot \Grad A } \dt \br 
    &\qquad + \int_0^\tau \intO{ \frac{\vr^1}{\Ov{\vr}} (\vU - \vu) \cdot \Grad A } \dt \br 
    &=\frac{1}{\eps} \int_0^\tau \intO{ \Big( \frac{\vr_\eps}{\Ov{\vr}} - 1 \Big) (\vu_\eps - \vU) \cdot \Grad A } \dt + \int_0^\tau \intO{ \frac{1}{\Ov{\vr}} ( \mathfrak{R} - \vr^1 ) \vU \cdot \Grad A } \dt \br 
    &\qquad + \int_0^\tau \intO{ \frac{\vr^1}{\Ov{\vr}} (\vU - \vu) \cdot \Grad A } \dt + O(\eps).
\end{align}

Hence,
\begin{align}
R &= - \int_0^\tau \intO{ \frac{\kappa(\Ov{\vt})}{c_p} \Big( \del_\vr s(\Ov{\vr}, \Ov{\vt}) (\mathfrak{R} - \vr^1) + \del_\vt s(\Ov{\vr}, \Ov{\vt}) (\mathfrak{T} - \vt^1) \Big) \Delta \vt^1 } \dt \br
&\qquad - \int_0^\tau \intO{ \frac{\Ov{\vr} \Ov{\vt} \alpha}{c_p} \Big( \del_\vr s(\Ov{\vr}, \Ov{\vt}) (\mathfrak{R} - \vr^1) + \del_\vt s(\Ov{\vr}, \Ov{\vt}) (\mathfrak{T} - \vt^1) \Big) \vU \cdot \Grad G } \dt \br
&\qquad - \int_0^\tau \intO{ \frac{\Ov{\vt} \alpha}{c_p} \Big( \del_\vr s(\Ov{\vr}, \Ov{\vt}) (\mathfrak{R} - \vr^1) + \del_\vt s(\Ov{\vr}, \Ov{\vt}) (\mathfrak{T} - \vt^1) \Big) \del_t \chi } \dt \br
    &\qquad + \int_0^\tau \intO{ \frac{\Ov{\vt} \alpha}{c_p} \Big( \del_\vr s(\Ov{\vr}, \Ov{\vt}) (\mathfrak{R} - \vr^1) + \del_\vt s(\Ov{\vr}, \Ov{\vt}) (\mathfrak{T} - \vt^1) \Big) ( \del_t A + \vU \cdot \Grad A ) } \br 
	&\qquad + \int_0^\tau \intO{ (\vr^1 - \mathfrak{R}) \vU \cdot \Grad G } \dt \br
	 	&\qquad + \int_0^\tau \intO{ \frac{1}{\Ov{\vr}} (\mathfrak{R} - \vr^1) (\del_t A + \vU \cdot \Grad A) } \dd t \br 
        &\qquad +\int_0^\tau \intO{ \frac{\vr^1}{\Ov{\vr}} (\vU - \vu) \cdot \Grad A } \dt \br 
        &\quad \quad + \frac{1}{\eps} \int_0^\tau \intO{ \Big( \frac{\vr_\eps}{\Ov{\vr}} - 1 \Big) ( \vu_\eps - \vU) \cdot \Grad A } \dt \br 
	&\qquad + \frac{1}{\eps^2} \int_0^\tau \intO{ (\vU \times (\vB_\eps - \vc H_\eps)) \cdot \Curl(\vB_\eps - \vc H_\eps) } \dd t \br 
	&\qquad + \frac{1}{\eps^2} \int_0^\tau \intO{ ((\vU - \vu_\eps) \times (\vB_\eps - \vc H_\eps)) \cdot \Curl \vc H_\eps } \dd t \br 
    &\qquad + \int_0^\tau \intO{ E_\eps\left( \vr_\eps, \vt_\eps, \vu_\eps, \vB_\eps \ \Big| r_\eps, \Theta_\eps, \vU, \vc H_\eps \right) } \dt + O(\eps). \label{rest-R1}
\end{align}

\paragraph{Step 5:} To handle the term containing $\del_t \chi$, we use $\int_\Omega \mathfrak{R} \dx = \int_\Omega \vr^1 \dx = 0$ to write using the Gibbs' equation
\begin{align*}
&- \int_0^\tau \intO{ \frac{\Ov{\vt} \alpha}{c_p} \Big( \del_\vr s(\Ov{\vr}, \Ov{\vt}) (\mathfrak{R} - \vr^1) + \del_\vt s(\Ov{\vr}, \Ov{\vt}) (\mathfrak{T} - \vt^1) \Big) \del_t \chi } \dt \br 
&\qquad = - \int_0^\tau \intO{ \frac{\Ov{\vt} \alpha}{c_p} \del_\vt s(\Ov{\vr}, \Ov{\vt}) (\mathfrak{T} - \vt^1) \del_t \chi } \dt \br 
&\qquad = - \int_0^\tau \intO{ \frac{\alpha \del_\vt e(\Ov{\vr}, \Ov{\vt})}{c_p} (\mathfrak{T} - \vt^1) \del_t \chi } \dt.
\end{align*}

For $\Delta_x \vt^1$ we find 
\begin{equation}\label{est-11}
\begin{split}
&- \int_0^\tau \intO{ \frac{\kappa(\Ov{\vt})}{c_p} \Big( \del_\vr s(\Ov{\vr}, \Ov{\vt}) (\mathfrak{R} - \vr^1) + \del_\vt s(\Ov{\vr}, \Ov{\vt}) (\mathfrak{T} - \vt^1) \Big) \Delta_x \vt^1 } \dt \\
&= - \int_0^\tau \intO{ \frac{\kappa(\Ov{\vt})}{c_p} \del_\vr s(\Ov{\vr}, \Ov{\vt}) \Big( (\mathfrak{R} - \vr^1) + \frac{\del_\vt p(\Ov{\vr}, \Ov{\vt})}{\del_\vr p(\Ov{\vr}, \Ov{\vt})} (\mathfrak{T} - \vt^1) + \frac{1}{\del_\vr p(\Ov{\vr}, \Ov{\vt})} (\mathfrak{A} - A) \Big) \Delta_x \vt^1 } \dt \\
&\qquad + \int_0^\tau \intO{ \frac{\kappa(\Ov{\vt})}{c_p} \Big( \del_\vr s(\Ov{\vr}, \Ov{\vt}) \frac{\del_\vt p(\Ov{\vr}, \Ov{\vt})}{\del_\vr p(\Ov{\vr}, \Ov{\vt})} - \del_\vt s(\Ov{\vr}, \Ov{\vt}) \Big) (\mathfrak{T} - \vt^1) \Delta_x \vt^1 } \dt \br 
&\qquad + \int_0^\tau \intO{ \frac{\kappa(\Ov{\vt})}{c_p} \frac{\del_\vr s(\Ov{\vr}, \Ov{\vt})}{\del_\vr p(\Ov{\vr}, \Ov{\vt})} (\mathfrak{A} - A) \Delta_x \vt^1 } \dt.
\end{split}
\end{equation}
Due to \eqref{BoussDiff}, the expression in brackets under the first integral is independent of $x$. Integrating the heat equation \eqref{target}$_3$ with respect to $x$, as $\Div \vU = 0$ and $\intO{A} = 0$, gives
\begin{align*}
\intO{ \Big(1 - \frac{\Ov{\vt} \alpha \del_\vt p(\Ov{\vr}, \Ov{\vt})}{\Ov{\vr} c_p} \Big) \frac{1}{\del_\vt p(\Ov{\vr}, \Ov{\vt})} \del_t \chi } = \Big(1 - \frac{\Ov{\vt} \alpha \del_\vt p(\Ov{\vr}, \Ov{\vt})}{\Ov{\vr} c_p} \Big) \int_\Omega \del_t \vt^1 \dd x = \int_\Omega \frac{\kappa(\Ov{\vt})}{\Ov{\vr} c_p} \Delta_x \vt^1 \dd x .
\end{align*}
Using the above we find using again $\intO{\mathfrak{R} - \vr^1} = 0 = \intO{\mathfrak{A} - A}$ that
\begin{align*}
&- \int_0^\tau \intO{ \frac{\kappa(\Ov{\vt})}{c_p} \del_\vr s(\Ov{\vr}, \Ov{\vt}) \Big( (\mathfrak{R} - \vr^1) + \frac{\del_\vt p(\Ov{\vr}, \Ov{\vt})}{\del_\vr p(\Ov{\vr}, \Ov{\vt})} (\mathfrak{T} - \vt^1) + \frac{1}{\del_\vr p(\Ov{\vr}, \Ov{\vt})} (\mathfrak{A} - A) \Big) \Delta_x \vt^1 } \dt \br 
&= - \int_0^\tau \intO{ \frac{\Ov{\vr} \del_\vr s(\Ov{\vr}, \Ov{\vt})}{\del_\vt p(\Ov{\vr}, \Ov{\vt})} \Big( 1 - \frac{\Ov{\vt} \alpha \del_\vt p(\Ov{\vr}, \Ov{\vt})}{\Ov{\vr} c_p} \Big) \Big( (\mathfrak{R} - \vr^1) + \frac{\del_\vt p(\Ov{\vr}, \Ov{\vt})}{\del_\vr p(\Ov{\vr}, \Ov{\vt})} (\mathfrak{T} - \vt^1) + \frac{1}{\del_\vr p(\Ov{\vr}, \Ov{\vt})} (\mathfrak{A} - A) \Big) \del_t \chi } \dd t \br 
&= - \int_0^\tau \intO{ \frac{\Ov{\vr} \del_\vr s(\Ov{\vr}, \Ov{\vt})}{\del_\vt p(\Ov{\vr}, \Ov{\vt})} \Big( 1 - \frac{\Ov{\vt} \alpha \del_\vt p(\Ov{\vr}, \Ov{\vt})}{\Ov{\vr} c_p} \Big) \frac{\del_\vt p(\Ov{\vr}, \Ov{\vt})}{\del_\vr p(\Ov{\vr}, \Ov{\vt})} (\mathfrak{T} - \vt^1) \del_t \chi } \dd t.
\end{align*}

Accordingly, collecting the coefficients of $(\mathfrak{T} - \vt^1) \del_t \chi$, we get in accordance with Gibbs' relation and the definition of $\alpha$ in \eqref{Coeff}
\begin{align*}
- \frac{\alpha \del_\vt e(\Ov{\vr}, \Ov{\vt})}{c_p} - \frac{\Ov{\vr} \del_\vr s(\Ov{\vr}, \Ov{\vt})}{\del_\vt p(\Ov{\vr}, \Ov{\vt})} \Big( 1 - \frac{\Ov{\vt} \alpha \del_\vt p(\Ov{\vr}, \Ov{\vt})}{\Ov{\vr} c_p} \Big) \frac{\del_\vt p(\Ov{\vr}, \Ov{\vt})}{\del_\vr p(\Ov{\vr}, \Ov{\vt})} = 0.
\end{align*}
Hence, all the terms containing $\del_t \chi$ add up to zero such that
\begin{equation}
\begin{split}
R &= \int_0^\tau \intO{ \frac{\kappa(\Ov{\vt})}{c_p} \Big( \del_\vr s(\Ov{\vr}, \Ov{\vt}) \frac{\del_\vt p(\Ov{\vr}, \Ov{\vt})}{\del_\vr p(\Ov{\vr}, \Ov{\vt})} - \del_\vt s(\Ov{\vr}, \Ov{\vt}) \Big) (\mathfrak{T} - \vt^1) \Delta_x \vt^1 } \dt \br 
&\qquad + \int_0^\tau \intO{ \frac{\kappa(\Ov{\vt})}{c_p} \frac{\del_\vr s(\Ov{\vr}, \Ov{\vt})}{\del_\vr p(\Ov{\vr}, \Ov{\vt})} (\mathfrak{A} - A) \Delta_x \vt^1 } \dt \br 
&\qquad - \int_0^\tau \intO{ \frac{\Ov{\vr} \Ov{\vt} \alpha}{c_p} \Big( \del_\vr s(\Ov{\vr}, \Ov{\vt}) (\mathfrak{R} - \vr^1) + \del_\vt s(\Ov{\vr}, \Ov{\vt}) (\mathfrak{T} - \vt^1) \Big) \vU \cdot \Grad G } \dt \br
    &\qquad + \int_0^\tau \intO{ \frac{\Ov{\vt} \alpha}{c_p} \Big( \del_\vr s(\Ov{\vr}, \Ov{\vt}) (\mathfrak{R} - \vr^1) + \del_\vt s(\Ov{\vr}, \Ov{\vt}) (\mathfrak{T} - \vt^1) \Big) ( \del_t A + \vU \cdot \Grad A ) } \br 
	&\qquad + \int_0^\tau \intO{ (\vr^1 - \mathfrak{R}) \vU \cdot \Grad G } \dt \br 
	 	&\qquad + \int_0^\tau \intO{ \frac{1}{\Ov{\vr}} (\mathfrak{R} - \vr^1) (\del_t A + \vU \cdot \Grad A) } \dd t \br 
        &\qquad + \int_0^\tau \intO{ \frac{\vr^1}{\Ov{\vr}} (\vU - \vu) \cdot \Grad A } \dt \br 
        &\quad \quad + \frac{1}{\eps} \int_0^\tau \intO{ \Big( \frac{\vr_\eps}{\Ov{\vr}} - 1 \Big) ( \vu_\eps - \vU) \cdot \Grad A } \dt \br 
	&\qquad + \frac{1}{\eps^2} \int_0^\tau \intO{ (\vU \times (\vB_\eps - \vc H_\eps)) \cdot \Curl(\vB_\eps - \vc H_\eps) } \dd t \br 
	&\qquad + \frac{1}{\eps^2} \int_0^\tau \intO{ ((\vU - \vu_\eps) \times (\vB_\eps - \vc H_\eps)) \cdot \Curl \vc H_\eps } \dd t \\
    &\qquad + \int_0^\tau \intO{ E_\eps\left( \vr_\eps, \vt_\eps, \vu_\eps, \vB_\eps \ \Big| r_\eps, \Theta_\eps, \vU, \vc H_\eps \right) } \dt + O(\eps). 
\end{split}
\end{equation}

\paragraph{Step 6:} Let us look on the $\eps$-independent terms containing $\vU \cdot \Grad G$. By definition of $c_p$ in \eqref{Coeff} and Boussinesq relation \eqref{BoussDiff} we see
\begin{align*}
& - \intO{ \frac{\Ov{\vr} \Ov{\vt} \alpha}{c_p} \Big( \del_\vr s(\Ov{\vr}, \Ov{\vt}) (\mathfrak{R} - \vr^1) + \del_\vt s(\Ov{\vr}, \Ov{\vt}) (\mathfrak{T} - \vt^1) \Big) \vU \cdot \Grad G } + \intO{ (\vr^1 - \mathfrak{R}) \Grad G \cdot \vU } \br 
&= \intO{ \frac{\Ov{\vr} \Ov{\vt} \alpha}{c_p} \Big( \del_\vr s(\Ov{\vr}, \Ov{\vt}) \Grad (\mathfrak{R} - \vr^1) + \del_\vt s(\Ov{\vr}, \Ov{\vt}) \Grad (\mathfrak{T} - \vt^1) \Big) \cdot \vU G } - \intO{ \Grad (\vr^1 - \mathfrak{R}) G \cdot \vU } \br 
&= - \intO{ \Big( \frac{\Ov{\vr} \Ov{\vt} \alpha}{c_p} \Big( \frac{\del_\vt p(\Ov{\vr}, \Ov{\vt})}{\Ov{\vr}^2} + \frac{\del_\vt e(\Ov{\vr}, \Ov{\vt}) \del_\vr p(\Ov{\vr}, \Ov{\vt})}{\Ov{\vt} \del_\vt p(\Ov{\vr}, \Ov{\vt})} \Big) \Grad (\mathfrak{R} - \vr^1) + \frac{\Ov{\vr} \alpha \del_\vt e(\Ov{\vr}, \Ov{\vt})}{c_p \del_\vt p(\Ov{\vr}, \Ov{\vt})} \Grad (\mathfrak{A} - A) \Big) \cdot \vU G } \br 
&\qquad - \intO{ \Grad (\vr^1 - \mathfrak{R}) G \cdot \vU } \br 
&= - \intO{ \frac{\Ov{\vr} \alpha \del_\vt e(\Ov{\vr}, \Ov{\vt})}{c_p \del_\vt p(\Ov{\vr}, \Ov{\vt})} \Grad (\mathfrak{A} - A) \cdot \vU G } = \intO{ \frac{\Ov{\vr} \alpha \del_\vt e(\Ov{\vr}, \Ov{\vt})}{c_p \del_\vt p(\Ov{\vr}, \Ov{\vt})} (\mathfrak{A} - A) \vU \cdot \Grad G } \br 
\end{align*}
and hence
\begin{align}
R &= \int_0^\tau \intO{ \frac{\kappa(\Ov{\vt})}{c_p} \Big( \del_\vr s(\Ov{\vr}, \Ov{\vt}) \frac{\del_\vt p(\Ov{\vr}, \Ov{\vt})}{\del_\vr p(\Ov{\vr}, \Ov{\vt})} - \del_\vt s(\Ov{\vr}, \Ov{\vt}) \Big) (\mathfrak{T} - \vt^1) \Delta_x \vt^1 } \dt \br
&\qquad + \int_0^\tau \intO{ \frac{\kappa(\Ov{\vt})}{c_p} \frac{\del_\vr s(\Ov{\vr}, \Ov{\vt})}{\del_\vr p(\Ov{\vr}, \Ov{\vt})} (\mathfrak{A} - A) \Delta_x \vt^1 } \dt \br 
&\qquad + \int_0^\tau \intO{ \frac{\Ov{\vr} \alpha \del_\vt e(\Ov{\vr}, \Ov{\vt})}{c_p \del_\vt p(\Ov{\vr}, \Ov{\vt})} (\mathfrak{A} - A) \vU \cdot \Grad G } \dt \br 
&\qquad + \int_0^\tau \intO{ \frac{\vr^1}{\Ov{\vr}} (\vU - \vu) \cdot \Grad A } \dt \br 
    &\qquad + \int_0^\tau \intO{ \frac{\Ov{\vt} \alpha}{c_p} \Big( \del_\vr s(\Ov{\vr}, \Ov{\vt}) (\mathfrak{R} - \vr^1) + \del_\vt s(\Ov{\vr}, \Ov{\vt}) (\mathfrak{T} - \vt^1) \Big) ( \del_t A + \vU \cdot \Grad A ) } \br 
	 	&\qquad + \int_0^\tau \intO{ \frac{1}{\Ov{\vr}} (\mathfrak{R} - \vr^1) (\del_t A + \vU \cdot \Grad A) } \dd t \br 
        &\quad \quad + \frac{1}{\eps} \int_0^\tau \intO{ \Big( \frac{\vr_\eps}{\Ov{\vr}} - 1 \Big) ( \vu_\eps - \vU) \cdot \Grad A } \dt \br 
	&\qquad + \frac{1}{\eps^2} \int_0^\tau \intO{ (\vU \times (\vB_\eps - \vc H_\eps)) \cdot \Curl(\vB_\eps - \vc H_\eps) } \dd t \br 
	&\qquad + \frac{1}{\eps^2} \int_0^\tau \intO{ ((\vU - \vu_\eps) \times (\vB_\eps - \vc H_\eps)) \cdot \Curl \vc H_\eps } \dd t \br 
    &\qquad + \int_0^\tau \intO{ E_\eps\left( \vr_\eps, \vt_\eps, \vu_\eps, \vB_\eps \ \Big| r_\eps, \Theta_\eps, \vU, \vc H_\eps \right) } \dt + O(\eps). \label{R1}
\end{align}

For the first term on the right hand side of the above we find by Gibbs' relation and the definition of $c_p$ in \eqref{Coeff} that 
\begin{align*}
\frac{\kappa(\Ov{\vt})}{c_p} \Big( \del_\vr s(\Ov{\vr}, \Ov{\vt}) \frac{\del_\vt p(\Ov{\vr}, \Ov{\vt})}{\del_\vr p(\Ov{\vr}, \Ov{\vt})} - \del_\vt s(\Ov{\vr}, \Ov{\vt}) \Big) = - \frac{\kappa(\Ov{\vt})}{\Ov{\vt}}.
\end{align*}

We will now focus on terms containing $A$. First,
\begin{align*}
    &\int_0^\tau \intO{ \frac{\vr^1}{\Ov{\vr}} (\vU - \vu) \cdot \Grad A } \dt + \frac{1}{\eps} \int_0^\tau \intO{ \Big( \frac{\vr_\eps}{\Ov{\vr}} - 1 \Big) (\vu_\eps - \vU) \cdot \Grad A } \dt \br 
    &= \int_0^\tau \intO{ \frac{1}{\Ov{\vr}} \Big( \vr^1 - \frac{\vr_\eps - \Ov{\vr}}{\eps} \Big) (\vU - \vu_\eps) \cdot \Grad A } \dt + O(\eps) \br 
    &= \int_0^\tau \intO{ \frac{1}{\Ov{\vr}} \frac{r_\eps - \vr_\eps}{\eps} (\vU - \vu_\eps) \cdot \Grad A } \dt + O(\eps) \br 
    &\lesssim \int_0^\tau \intO{ E_\eps\left( \vr_\eps, \vt_\eps, \vu_\eps, \vB_\eps \ \Big| r_\eps, \Theta_\eps, \vU, \vc H_\eps \right) } \dt + O(\eps).
\end{align*}

The fifth and sixth term in \eqref{R1} can be transformed into
\begin{align*}
    &\int_0^\tau \intO{ \frac{\alpha}{c_p} \Big( - \frac{\Ov{\vt}}{\Ov{\vr}^2} \del_\vt p(\Ov{\vr}, \Ov{\vt}) (\mathfrak{R} - \vr^1) + \del_\vt e(\Ov{\vr}, \Ov{\vt}) (\mathfrak{T} - \vt^1) \Big) ( \del_t A + \vU \cdot \Grad A ) } \dt \br 
    &\qquad + \int_0^\tau \intO{ \frac{1}{\Ov{\vr}} (\mathfrak{R} - \vr^1) (\del_t A + \vU \cdot \Grad A) } \dd t \br 
    &= \int_0^\tau \intO{ \Big( \Big( \frac{1}{\Ov{\vr}} - \frac{\alpha}{c_p} \frac{\Ov{\vt}}{\Ov{\vr}^2} \del_\vt p(\Ov{\vr}, \Ov{\vt}) \Big) (\mathfrak{R} - \vr^1) + \frac{\alpha \del_\vt e(\Ov{\vr}, \Ov{\vt})}{c_p} (\mathfrak{T} - \vt^1) \Big) ( \del_t A + \vU \cdot \Grad A ) } \dt \br 
    &= \int_0^\tau \intO{ \frac{\del_\vt e(\Ov{\vr}, \Ov{\vt})}{\Ov{\vr} c_p} \Big( (\mathfrak{R} - \vr^1) + \frac{\del_\vt p(\Ov{\vr}, \Ov{\vt})}{\del_\vr p(\Ov{\vr}, \Ov{\vt})} (\mathfrak{T} - \vt^1) \Big) ( \del_t A + \vU \cdot \Grad A ) } \dt \br 
    &= \int_0^\tau \intO{ \frac{\del_\vt e(\Ov{\vr}, \Ov{\vt})}{\Ov{\vr} c_p} \Big( (\mathfrak{R} - \vr^1) + \frac{\del_\vt p(\Ov{\vr}, \Ov{\vt})}{\del_\vr p(\Ov{\vr}, \Ov{\vt})} (\mathfrak{T} - \vt^1) + \frac{1}{\del_\vr p(\Ov{\vr}, \Ov{\vt})} (\mathfrak{A} - A) \Big) ( \del_t A + \vU \cdot \Grad A ) } \dt \br 
    &\qquad - \int_0^\tau \intO{ \frac{\del_\vt e(\Ov{\vr}, \Ov{\vt})}{\Ov{\vr} c_p \del_\vr p(\Ov{\vr}, \Ov{\vt})} (\mathfrak{A} - A) ( \del_t A + \vU \cdot \Grad A ) } \dt \br 
    &= - \int_0^\tau \intO{ \frac{\del_\vt e(\Ov{\vr}, \Ov{\vt})}{\Ov{\vr} c_p \del_\vr p(\Ov{\vr}, \Ov{\vt})} (\mathfrak{A} - A) ( \del_t A + \vU \cdot \Grad A ) } \dt
\end{align*}
since the last expression in brackets is independent of $x$, $\Div \vU = 0$, and $\intO{\del_t A} = \del_t \intO{A} = 0$ since $A$ is mean-free. We will show that the terms containing $(\mathfrak{A} - A)$ vanish identically. Indeed, as $A = \Ov{\vB} \cdot \vB^1$ satisfies induction equation, we may write this as (using $\alpha = \Ov{\vr}^{-1} \del_\vt p(\Ov{\vr}, \Ov{\vt}) / \del_\vr p(\Ov{\vr}, \Ov{\vt})$)
\begin{align*}
    &- \int_0^\tau \intO{ \frac{1}{\Ov{\vr}^2 c_p \del_\vr p(\Ov{\vr}, \Ov{\vt})} \Big( \Ov{\vr} \del_\vt e(\Ov{\vr}, \Ov{\vt}) \zeta(\Ov{\vt}) \Delta_x A + \kappa(\Ov{\vt}) \del_\vt p(\Ov{\vr}, \Ov{\vt}) \Delta_x \vt^1 - \Ov{\vr}^2 \del_\vt e(\Ov{\vr}, \Ov{\vt}) \vU \cdot \Grad G \Big) (\mathfrak{A} - A) } \dt.
\end{align*}
Going further and multiplying in \eqref{target} the induction equation by $\Ov{\vr} \del_\vt e(\Ov{\vr}, \Ov{\vt})$ and the heat equation by $\del_\vt p(\Ov{\vr}, \Ov{\vt})$, summing up, multiplying by $(\mathfrak{A} - A)$ and integrating, and taking also into account that the occurring integral of $(\mathfrak{A} - A)\del_t \chi$ will vanish since $\intO{\mathfrak{A} - A} = 0$, we get together with the definition of $c_p = \del_\vt e(\Ov{\vr}, \Ov{\vt}) + \Ov{\vt} \Ov{\vr}^{-1} \alpha \del_\vt p(\Ov{\vr}, \Ov{\vt})$
\begin{align}
    &\int_0^\tau \intO{ (\Ov{\vr} \del_\vt e(\Ov{\vr}, \Ov{\vt}) \zeta(\Ov{\vt}) \Delta_x A + \kappa(\Ov{\vt}) \del_\vt p(\Ov{\vr}, \Ov{\vt}) \Delta_x \vt^1 - \Ov{\vr}^2 \del_\vt e(\Ov{\vr}, \Ov{\vt}) \vU \cdot \Grad G ) (\mathfrak{A} - A)} \dt \br 
    &= \int_0^\tau \intO{ [ (\Ov{\vr} \del_\vt e(\Ov{\vr}, \Ov{\vt}) + \del_\vt p(\Ov{\vr}, \Ov{\vt}) \Ov{\vt} \alpha) (\del_t A + \vU \cdot \Grad A) + \del_\vt p(\Ov{\vr}, \Ov{\vt}) \Ov{\vr} c_p (\del_t \vt^1 + \vU \cdot \Grad \vt^1)] (\mathfrak{A} - A) } \dt \br 
    &\qquad - \int_0^\tau \intO{ (\del_\vt p(\Ov{\vr}, \Ov{\vt}) \Ov{\vr} \Ov{\vt} \alpha + \Ov{\vr}^2 \del_\vt e(\Ov{\vr}, \Ov{\vt})) \vU \cdot \Grad G ) (\mathfrak{A} - A) } \dt \br 
    &= \int_0^\tau \intO{ \Ov{\vr} c_p [ (\del_t A + \vU \cdot \Grad A) + \del_\vt p(\Ov{\vr}, \Ov{\vt}) (\del_t \vt^1 + \vU \cdot \Grad \vt^1) - \Ov{\vr} \vU \cdot \Grad G ] (\mathfrak{A} - A) } \dt \br 
    &= \int_0^\tau \intO{ \Ov{\vr} c_p [ \del_t ( A + \del_\vt p(\Ov{\vr}, \Ov{\vt}) \vt^1 ) + \vU \cdot \Grad ( A + \del_\vt p(\Ov{\vr}, \Ov{\vt}) \vt^1 - \Ov{\vr} G) ] (\mathfrak{A} - A) } \dt.
\end{align}
Using further Boussinesq relation \eqref{OB}, we find
\begin{align}
    &\int_0^\tau \intO{ \Ov{\vr} c_p [ \del_t ( A + \del_\vt p(\Ov{\vr}, \Ov{\vt}) \vt^1 ) + \vU \cdot \Grad ( A + \del_\vt p(\Ov{\vr}, \Ov{\vt}) \vt^1 - \Ov{\vr} G) ] (\mathfrak{A} - A) } \dt \br 
    &= \int_0^\tau \intO{ \Ov{\vr} c_p [ \del_t ( \chi - \del_\vr p(\Ov{\vr}, \Ov{\vt}) \vr^1 ) - \del_\vr p(\Ov{\vr}, \Ov{\vt}) \vU \cdot \Grad \vr^1 ] (\mathfrak{A} - A) } \dt \br 
    &= - \int_0^\tau \intO{ \Ov{\vr} c_p \del_\vr p(\Ov{\vr}, \Ov{\vt}) [ \del_t \vr^1 + \vU \cdot \Grad \vr^1 ] (\mathfrak{A} - A) } \dt \br 
    &= - \int_0^\tau \intO{ \Ov{\vr} c_p \del_\vr p(\Ov{\vr}, \Ov{\vt}) [ \del_t \vr^1 + \Div(\vr^1 \vU) ] (\mathfrak{A} - A) } \dt = 0,
\end{align}
the last line coming from the fact that $\vr^1$ fulfills the continuity equation \eqref{target}$_2$. In turn, the remainder collapses to
\begin{align*}
R &= - \int_0^\tau \intO{ \frac{\kappa(\Ov{\vt})}{\Ov{\vt}} (\mathfrak{T} - \vt^1) \Delta_x \vt^1 } \dt \br 
	&\qquad + \frac{1}{\eps^2} \int_0^\tau \intO{ (\vU \times (\vB_\eps - \vc H_\eps)) \cdot \Curl(\vB_\eps - \vc H_\eps) } \dd t \br 
	&\qquad + \frac{1}{\eps^2} \int_0^\tau \intO{ ((\vU - \vu_\eps) \times (\vB_\eps - \vc H_\eps)) \cdot \Curl \vc H_\eps } \dd t \\
    &\qquad + \int_0^\tau \intO{ E_\eps\left( \vr_\eps, \vt_\eps, \vu_\eps, \vB_\eps \ \Big| r_\eps, \Theta_\eps, \vU, \vc H_\eps \right) } \dt + O(\eps).
\end{align*}

Consequently, we may write the relative energy inequality \eqref{ei2} in the form
\begin{align}
& \left[ \intO{ E_\eps\left( \vr_\eps, \vt_\eps, \vu_\eps, \vB_\eps \ \Big| r_\eps, \Theta_\eps, \vU, \vc H_\eps \right) } \right]_{t = 0}^{t = \tau} \br 
&+ \int_0^\tau \intO{ 	\frac{\Theta_\eps}{\vt_\eps} \left( \S(\vt_\eps, \Grad \vu_\eps) : \Grad \vu_\eps + \frac{1}{\eps^2} \frac{\kappa(\vt_\eps) |\Grad \vt_\eps|^2}{\vt_\eps} + \frac{1}{\eps^2} \zeta(\vt_\eps) |\Curl \vB_\eps|^2 \right) } \dt \br 
	&- \int_0^\tau \intO{ \S(\Ov{\vt}, \Grad \vU) : \Grad (\vu_\eps - \vU) } \dt - \int_0^\tau \intO{ \S(\vt_\eps, \nabla \vu_\eps) : \Grad \vU } \dt \br 
	&- \frac{1}{\eps^2} \int_0^\tau \intO{ \frac{\kappa(\vt_\eps) }{ \vt_\eps} \Grad \vt_\eps \cdot \Grad \Theta_\eps } - \int_0^\tau \intO{ \frac{\kappa(\Ov{\vt})}{\Ov{\vt}} \Grad (\mathfrak{T} - \vt^1) \cdot \Grad \vt^1 } \dt \br 
	&- \frac{1}{\eps^2} \int_0^\tau \intO{ \zeta(\Ov{\vt}) \Curl (\vB_\eps - \vc H_\eps) \cdot \Curl \vc H_\eps + \zeta(\vt_\eps) \Curl \vB_\eps \cdot \Curl \vc H_\eps     } \dt \br 
	&\leq \frac{1}{\eps^2} \int_0^\tau \intO{ (\vU \times (\vB_\eps - \vc H_\eps)) \cdot \Curl(\vB_\eps - \vc H_\eps) } \dt \br 
	&\qquad + \frac{1}{\eps^2} \int_0^\tau \intO{ ((\vU - \vu_\eps) \times (\vB_\eps - \vc H_\eps)) \cdot \Curl \vc H_\eps } \dt \br  
	&\qquad + \int_0^\tau \intO{ E_\eps\left( \vr_\eps, \vt_\eps, \vu_\eps, \vB_\eps \ \Big| r_\eps, \Theta_\eps, \vU, \vc H_\eps \right) } \dt + O(\eps). \label{est-12}
\end{align}

\paragraph{Step 7:} We want now to manipulate the left hand site of the foregoing inequality to get non-negative diffusive terms. To this end, for the temperature, using strong convergence of $\vt_\eps \to \Ov{\vt}$ and \eqref{w-conv_T}, we may write
\begin{align}\label{est-14}
&\int_0^\tau \intO{ \frac{1}{\eps^2} \frac{\Theta_\eps}{\vt_\eps} \frac{\kappa(\vt_\eps)}{\vt_\eps} \Grad \vt_\eps \cdot \Grad \vt_\eps - \frac{1}{\eps^2} \frac{\kappa(\vt_\eps)}{\vt_\eps} \Grad \vt_\eps \cdot \Grad \Theta_\eps - \frac{\kappa(\Ov{\vt})}{\Ov{\vt}} \Grad (\mathfrak{T} - \vt^1) \cdot \Grad \vt^1 } \dt \br 
&= \int_0^\tau \intO{ \frac{1}{\eps^2} \frac{\Theta_\eps}{\vt_\eps} \frac{\kappa(\vt_\eps)}{\vt_\eps} \Grad \vt_\eps \cdot \Grad \vt_\eps - \frac{1}{\eps^2} \frac{\kappa(\vt_\eps)}{\vt_\eps} \Grad \vt_\eps \cdot \Grad \Theta_\eps - \frac{1}{\eps^2} \frac{\kappa(\Ov{\vt})}{\Ov{\vt}} \Grad \vt_\eps \cdot \Grad \Theta_\eps } \dt \br 
&\quad + \frac{1}{\eps^2} \int_0^\tau \intO{ \frac{\kappa(\Ov{\vt})}{\Ov{\vt}} \Grad \Theta_\eps \cdot \Grad \Theta_\eps } \dt + O(\eps) \br 
&= \frac{1}{\eps^2} \int_0^\tau \intO{ \Big( \frac{\Theta_\eps}{\vt_\eps} - 1 \Big) \frac{\kappa(\vt_\eps)}{\vt_\eps} \Grad \vt_\eps \cdot \Grad \vt_\eps + \frac{1}{\eps^2} \Big( \frac{\vt_\eps}{\Theta_\eps} - 1 \Big) \frac{\kappa(\Theta_\eps)}{\Theta_\eps} \Grad \Theta_\eps \cdot \Grad \Theta_\eps } \dt \br 
&\quad + \frac{1}{\eps^2} \int_0^\tau \intO{ \Big( \frac{\kappa(\Theta_\eps)}{\Theta_\eps} \Grad \Theta_\eps - \frac{\kappa(\vt_\eps)}{\vt_\eps} \Grad \vt_\eps \Big) \cdot (\Grad \Theta_\eps - \Grad \vt_\eps) } \dt \br 
&\quad - \frac{1}{\eps^2} \int_0^\tau \intO{ \Big( \frac{\vt_\eps}{\Theta_\eps} - 1 \Big) \frac{\kappa(\Theta_\eps)}{\Theta_\eps} \Grad \Theta_\eps \cdot \Grad \Theta_\eps } \dt - \frac{1}{\eps^2} \int_0^\tau \intO{ \frac{\kappa(\Theta_\eps)}{\Theta_\eps} \Grad \Theta_\eps \cdot \Grad (\Theta_\eps - \vt_\eps) } \dt \br 
&\quad - \frac{1}{\eps^2} \int_0^\tau \intO{ \frac{\kappa(\Ov{\vt})}{\Ov{\vt}} \Grad \vt_\eps \cdot \Grad \Theta_\eps } \dt + \frac{1}{\eps^2} \int_0^\tau \intO{ \frac{\kappa(\Ov{\vt})}{\Ov{\vt}} \Grad \Theta_\eps \cdot \Grad \Theta_\eps } \dt + O(\eps) \br 
&= \frac{1}{\eps^2} \int_0^\tau \intO{ \Big( \frac{\Theta_\eps}{\vt_\eps} - 1 \Big) \frac{\kappa(\vt_\eps)}{\vt_\eps} \Grad \vt_\eps \cdot \Grad \vt_\eps } \dt + \frac{1}{\eps^2} \int_0^\tau \intO{ \Big( \frac{\vt_\eps}{\Theta_\eps} - 1 \Big) \frac{\kappa(\Theta_\eps)}{\Theta_\eps} \Grad \Theta_\eps \cdot \Grad \Theta_\eps } \dt \br 
&\quad + \frac{1}{\eps^2} \int_0^\tau \intO{ \Big( \frac{\kappa(\Theta_\eps)}{\Theta_\eps} \Grad \Theta_\eps - \frac{\kappa(\vt_\eps)}{\vt_\eps} \Grad \vt_\eps \Big) \cdot (\Grad \Theta_\eps - \Grad \vt_\eps) } \dt \br 
&\quad - \frac{1}{\eps^2} \int_0^\tau \intO{ \Big( \frac{\vt_\eps}{\Theta_\eps} - 1 \Big) \frac{\kappa(\Theta_\eps)}{\Theta_\eps} \Grad \Theta_\eps \cdot \Grad \Theta_\eps } \dt \br 
&\qquad + \frac{1}{\eps^2} \int_0^\tau \intO{ \Big( \frac{\kappa(\Theta_\eps)}{\Theta_\eps} - \frac{\kappa(\Ov{\vt})}{\Ov{\vt}} \Big) ( \Grad \vt_\eps - \Grad \Theta_\eps) \cdot \Grad \Theta_\eps } \dt + O(\eps)
\br  
&= \frac{1}{\eps^2} \int_0^\tau \intO{ \Big( \frac{\Theta_\eps}{\vt_\eps} - 1 \Big) \frac{\kappa(\vt_\eps)}{\vt_\eps} \Grad \vt_\eps \cdot \Grad \vt_\eps + \frac{1}{\eps^2} \Big( \frac{\vt_\eps}{\Theta_\eps} - 1 \Big) \frac{\kappa(\Theta_\eps)}{\Theta_\eps} \Grad \Theta_\eps \cdot \Grad \Theta_\eps } \dt \br 
&\quad + \frac{1}{\eps^2} \int_0^\tau \intO{ \Big( \frac{\kappa(\Theta_\eps)}{\Theta_\eps} \Grad \Theta_\eps - \frac{\kappa(\vt_\eps)}{\vt_\eps} \Grad \vt_\eps \Big) \cdot (\Grad \Theta_\eps - \Grad \vt_\eps) } \dt + O(\eps) .
\end{align}
Next, the viscous stress becomes
\begin{align}\label{est-15}
&\int_0^\tau \intO{ \frac{\Theta_\eps}{\vt_\eps} \S(\vt_\eps, \Grad \vu_\eps) : \Grad \vu_\eps - \S(\Ov{\vt}, \Grad \vU) : \Grad(\vu_\eps - \vU) - \S(\vt_\eps, \Grad \vu_\eps) : \Grad \vU } \dt \br 
&\quad - \int_0^\tau \intO{ \S(\Ov{\vt}, \Grad \vU) : \Grad (\vu_\eps - \vU) } \dt + \int_0^\tau \intO{ \S(\Theta_\eps, \Grad \vU) : \Grad (\vu_\eps - \vU) } \dt \br 
&= \int_0^\tau \intO{ \Big( \frac{\Theta_\eps}{\vt_\eps} - 1 \Big) \S(\vt_\eps, \Grad \vu_\eps) : \Grad \vu_\eps + \Big( \frac{\vt_\eps}{\Theta_\eps} - 1 \Big) \S(\Theta_\eps, \Grad \vU) : \Grad \vU } \dt \br 
&\quad + \int_0^\tau \intO{ (\S(\vt_\eps, \Grad \vu_\eps) - \S(\Theta_\eps, \Grad \vU)) : \Grad (\vu_\eps - \vU) } \dt \br 
&\quad - \int_0^\tau \intO{ \Big( \frac{\vt_\eps}{\Theta_\eps} - 1 \Big) \S(\Theta_\eps, \Grad \vU) : \Grad \vU } \dt \br 
&\qquad + \int_0^\tau \intO{ ( \S(\Theta_\eps, \Grad \vU) - \S(\Ov{\vt}, \Grad \vU) ) : \Grad (\vu_\eps - \vU) } \dt \br 
&= \int_0^\tau \intO{ \Big( \frac{\Theta_\eps}{\vt_\eps} - 1 \Big) \S(\vt_\eps, \Grad \vu_\eps) : \Grad \vu_\eps + \Big( \frac{\vt_\eps}{\Theta_\eps} - 1 \Big) \S(\Theta_\eps, \Grad \vU) : \Grad \vU } \dt \br 
&\quad + \int_0^\tau \intO{ (\S(\vt_\eps, \Grad \vu_\eps) - \S(\Theta_\eps, \Grad \vU)) : \Grad (\vu_\eps - \vU) } \dt + O(\eps).
\end{align}

Finally, the magnetic field becomes
\begin{align}\label{est-16}
&\frac{1}{\eps^2} \int_0^\tau \intO{ \frac{\Theta_\eps}{\vt_\eps} \zeta(\vt_\eps) \Curl \vB_\eps \cdot \Curl \vB_\eps - \frac{1}{\eps^2} \zeta(\Ov{\vt}) \Curl(\vB_\eps - \vc H_\eps) \cdot \Curl \vc H_\eps } \dt \br 
&\quad - \frac{1}{\eps^2} \int_0^\tau \intO{ \zeta(\vt_\eps) \Curl \vB_\eps \cdot \Curl \vc H_\eps } \dt \br 
&= \frac{1}{\eps^2} \int_0^\tau \intO{ \Big( \frac{\Theta_\eps}{\vt_\eps} - 1 \Big) \zeta(\vt_\eps) \Curl \vB_\eps \cdot \Curl \vB_\eps } \dt \br 
&\qquad + \frac{1}{\eps^2} \int_0^\tau \intO{ \Big( \frac{\vt_\eps}{\Theta_\eps} - 1 \Big) \zeta(\Theta_\eps) \Curl \vc H_\eps \cdot \Curl \vc H_\eps } \dt \br 
&\quad + \frac{1}{\eps^2} \int_0^\tau \intO{ (\zeta(\vt_\eps) \Curl \vB_\eps - \zeta(\Theta_\eps) \Curl \vc H_\eps) \cdot \Curl (\vB_\eps - \vc H_\eps) } \dt \br 
&\quad - \frac{1}{\eps^2} \int_0^\tau \intO{ \Big( \frac{\vt_\eps}{\Theta_\eps} - 1 \Big) \zeta(\Theta_\eps) \Curl \vc H_\eps \cdot \Curl \vc H_\eps } \dt \br 
&\qquad + \frac{1}{\eps^2} \int_0^\tau \intO{ ( \zeta(\Theta_\eps) - \zeta(\Ov{\vt}) ) \Curl \vc H_\eps \cdot \Curl(\vB_\eps - \vc H_\eps) } \dt \br 
&= \frac{1}{\eps^2} \int_0^\tau \intO{ \Big( \frac{\Theta_\eps}{\vt_\eps} - 1 \Big) \zeta(\vt_\eps) \Curl \vB_\eps \cdot \Curl \vB_\eps } \dt \br 
&\qquad + \int_0^\tau \intO{ \frac{1}{\eps^2} \Big( \frac{\vt_\eps}{\Theta_\eps} - 1 \Big) \zeta(\Theta_\eps) \Curl \vc H_\eps \cdot \Curl \vc H_\eps } \dt \br 
&\quad + \frac{1}{\eps^2} \int_0^\tau \intO{ (\zeta(\vt_\eps) \Curl \vB_\eps - \zeta(\Theta_\eps) \Curl \vc H_\eps) \cdot \Curl (\vB_\eps - \vc H_\eps) } \dt + O(\eps).
\end{align}

Following now \cite[Section~4.2]{FeireislNovotny2022}, we have
\begin{align}
&\|\Grad (\vu_\eps - \vU) \|_{L^2((0,T) \times \Omega)}^2 + \frac{1}{\eps^2} \|\Grad (\vt_\eps - \Theta_\eps)\|_{L^2((0,T) \times \Omega)}^2 + \frac{1}{\eps^2} \|\Curl(\vB_\eps - \vc H_\eps)\|_{L^2((0,T) \times \Omega)}^2 \br 
&\lesssim \int_0^\tau \intO{ \Big( \frac{\Theta_\eps}{\vt_\eps} - 1 \Big) \S(\vt_\eps, \Grad \vu_\eps) : \Grad \vu_\eps + \Big( \frac{\vt_\eps}{\Theta_\eps} - 1 \Big) \S(\Theta_\eps, \Grad \vU) : \Grad \vU } \dt \br 
& + \int_0^\tau \intO{ (\S(\vt_\eps, \Grad \vu_\eps) - \S(\Theta_\eps, \Grad \vU)) : \Grad (\vu_\eps - \vU) } \dt \br 
&+ \frac{1}{\eps^2} \int_0^\tau \intO{ \Big( \frac{\Theta_\eps}{\vt_\eps} - 1 \Big) \frac{\kappa(\vt_\eps)}{\vt_\eps} \Grad \vt_\eps \cdot \Grad \vt_\eps + \Big( \frac{\vt_\eps}{\Theta_\eps} - 1 \Big) \frac{\kappa(\Theta_\eps)}{\Theta_\eps} \Grad \Theta_\eps \cdot \Grad \Theta_\eps } \dt \br 
&+ \frac{1}{\eps^2} \int_0^\tau \intO{ \Big( \frac{\kappa(\Theta_\eps)}{\Theta_\eps} \Grad \Theta_\eps - \frac{\kappa(\vt_\eps)}{\vt_\eps} \Grad \vt_\eps \Big) \cdot (\Grad \Theta_\eps - \Grad \vt_\eps) } \dt \br 
&+ \frac{1}{\eps^2} \int_0^\tau \intO{ \Big( \frac{\Theta_\eps}{\vt_\eps} - 1 \Big) \zeta(\vt_\eps) \Curl \vB_\eps \cdot \Curl \vB_\eps + \Big( \frac{\vt_\eps}{\Theta_\eps} - 1 \Big) \zeta(\Theta_\eps) \Curl \vc H_\eps \cdot \Curl \vc H_\eps } \dt \br 
&+ \frac{1}{\eps^2} \int_0^\tau \intO{ (\zeta(\vt_\eps) \Curl \vB_\eps - \zeta(\Theta_\eps) \Curl \vc H_\eps) \cdot \Curl (\vB_\eps - \vc H_\eps) } \dt \br 
	&+ \int_0^\tau \intO{E_\eps\left( \vr_\eps, \vt_\eps, \vu_\eps, \vB_\eps \ \Big| r_\eps, \Theta_\eps, \vU, \vc H_\eps \right)} \dt .\label{est-13}
\end{align}

Summarizing  \eqref{est-12}, \eqref{est-14}, \eqref{est-15}, \eqref{est-16}, and \eqref{est-13} we can find the following  bound
\begin{align}\label{ei3}
& \left[ \intO{ E_\eps\left( \vr_\eps, \vt_\eps, \vu_\eps, \vB_\eps \ \Big| r_\eps, \Theta_\eps, \vU, \vc H_\eps \right) } \right]_{t = 0}^{t = \tau}
	\br 
	& + \|\Grad (\vu_\eps - \vU) \|_{L^2((0,T) \times \Omega)}^2 + \frac{1}{\eps^2} \|\Grad (\vt_\eps - \Theta_\eps)\|_{L^2((0,T) \times \Omega)}^2 + \frac{1}{\eps^2} \|\Curl(\vB_\eps - \vc H_\eps)\|_{L^2((0,T) \times \Omega)}^2 \br 
	&\quad \lesssim \frac{1}{\eps^2} \int_0^\tau \intO{ (\vU \times (\vB_\eps - \vc H_\eps)) \cdot \Curl(\vB_\eps - \vc H_\eps) } \dt \br 
	&\qquad + \frac{1}{\eps^2} \int_0^\tau \intO{ ((\vU - \vu_\eps) \times (\vB_\eps - \vc H_\eps)) \cdot \Curl \vc H_\eps } \dt \br 
	&\qquad + \int_0^\tau \intO{ E_\eps\left( \vr_\eps, \vt_\eps, \vu_\eps, \vB_\eps \ \Big| r_\eps, \Theta_\eps, \vU, \vc H_\eps \right) } \dt + O(\eps).
\end{align}

\paragraph{Step 8:} As for the last two magnetic integrals in \eqref{ei3}, we get (using \eqref{ess_rel})
\begin{align*}
&\frac{1}{\eps^2} \int_0^\tau \intO{ (\vU \times (\vB_\eps - \vc H_\eps)) \cdot \Curl(\vB_\eps - \vc H_\eps) } \dt \br 
	&\quad + \frac{1}{\eps^2} \int_0^\tau \intO{ ((\vU - \vu_\eps) \times (\vB_\eps - \vc H_\eps)) \cdot \Curl \vc H_\eps } \dt \br 
	&\lesssim \int_0^\tau \frac{1}{\eps} \|\vB_\eps - \vc H_\eps\|_{L^2(\Omega)} \cdot \frac{1}{\eps} \|\Curl(\vB_\eps - \vc H_\eps)\|_{L^2(\Omega)} \dd t \br 
	&\quad + \int_0^\tau \|\vu_\eps - \vU\|_{L^2(\Omega)} \cdot \frac{1}{\eps} \|\vB_\eps - \vc H_\eps\|_{L^2(\Omega)} \dt \br 
	&\lesssim \int_0^\tau \|\vu_\eps - \vU\|_{L^2(\Omega)}^2 + \frac{1}{\eps^2} \|\vB_\eps - \vc H_\eps\|_{L^2(\Omega)}^2 \dt + \delta \frac{1}{\eps^2} \|\Curl(\vB_\eps - \vc H_\eps)\|_{L^2((0,T) \times \Omega)}^2 \br 
	&\lesssim \int_0^\tau \intO{ E_\eps\left( \vr_\eps, \vt_\eps, \vu_\eps, \vB_\eps \ \Big| r_\eps, \Theta_\eps, \vU, \vc H_\eps \right) } \dt + \delta \frac{1}{\eps^2} \|\Curl(\vB_\eps - \vc H_\eps)\|_{L^2((0,T) \times \Omega)}^2 
\end{align*}
 for $\delta>0$ sufficiently small such that this term can be absorbed by the left hand site of \eqref{ei3}. Consequently, by \eqref{ei3} and the above, our final estimate reads
\begin{align*}
& \left[ \intO{ E_\eps\left( \vr_\eps, \vt_\eps, \vu_\eps, \vB_\eps \ \Big| r_\eps, \Theta_\eps, \vU, \vc H_\eps \right) } \right]_{t = 0}^{t = \tau}
	\br 
	& + \|\Grad (\vu_\eps - \vU) \|_{L^2((0,T) \times \Omega)}^2 + \frac{1}{\eps^2} \|\Grad (\vt_\eps - \Theta_\eps)\|_{L^2((0,T) \times \Omega)}^2 + \frac{1}{\eps^2} \|\Curl(\vB_\eps - \vc H_\eps)\|_{L^2((0,T) \times \Omega)}^2 \br 
	&\quad \lesssim \int_0^\tau \intO{ E_\eps\left( \vr_\eps, \vt_\eps, \vu_\eps, \vB_\eps \ \Big| r_\eps, \Theta_\eps, \vU, \vc H_\eps \right) } \dt + O(\eps).
\end{align*}

An application of Gr\"onwall's inequality now yields the desired conclusion \eqref{en_to_0} if \eqref{in_en_to_0} holds. The proof of Theorem~\ref{thm:main} is finished.






\section{Acknowledgements}
{\it The authors wish to thank Krzysztof Mizerski for many fruitful discussions.}

\paragraph{Conflict of interest:}
	The authors declare no conflict of interest in this paper.
	
\paragraph{Data availability:}
	There are no data available.




\bibliographystyle{amsplain}
\bibliography{Lit}

\end{document}